\documentclass[numbib]{imaiai}

\usepackage{amsmath,amsthm,amsfonts,amscd,amssymb,amsbsy,epsf}
\usepackage{verbatim}
\usepackage{graphics, subfigure}
\usepackage{bm}
\usepackage{enumerate}
\usepackage{bbm, dsfont}
\usepackage{array, contour}

\let\jmath=\undefined
\DeclareSymbolFont{cmletters}{OML}{cmm}{m}{it}
\DeclareMathSymbol{\jmath}{\mathord}{cmletters}{"7C}
\def\jversor{\contour[2]{black}{$\widehat{\jmath}$}}

\def\u{\mathbf{u} }
\def\v{\mathbf{v} }
\def\b{\mathbf{b} }
\def\a{\mathbf{a} }

\def\iversor{\bm{\widehat{\imath} } }
\def\kversor{\bm {\widehat{k} } }
\def\eversor{\bm {\widehat{e} } }

\def\n{\bm{\widehat{n} } }
\def\fversor{\bm {\widehat{f} } }
\def\zero{\bm { 0 } }

\def\A{\mathcal{A} }

\def\E{\mathcal{E} }
\def\F{\mathcal{F} }

\def\triang{\mathcal{T}}
\def\control{\mathcal{P}}
\def\D{\mathcal{D}}
\def\Hp{\mathcal{H}}
\def\Kz{K^{\zeta}}
\def\Lz{L^{\zeta}}

\def\pk{\mathbf{p}^{K}}
\def\pl{\mathbf{p}^{L}}

\def\qs{\mathbf{q}^{\sigma}}

\def\nukl{\nu^{\,K}_{\ell}}
\def\nusigk{\nu^{\,K}_{\sigma}}

\def\nutauk{\nu^{\,K}_{\tau}}

\def\ugamma{\mathbf{u}^{\Gamma}}
\def\vgamma{\mathbf{v}^{\Gamma}}
\def\zmax{z_{\max}}

\def\prob{\mathbbm{P}}
\def\entro{\mathbbm{H}}
\def\echoice{\mathbbm{h}}
\def\ind{\bm {\mathbbm{1} } }

\def\R{\bm{\mathbbm{R} } }
\def\C{\bm{\mathbbm{C} } }
\def\defining{\overset{\textbf{def}}=}

\begin{document}

\title{On the Construction of Geometric Parameters for \\
Preferential Fluid Flow Information in Fissured Media}

\shorttitle{Geometric Parameters and Preferential Flow Information} 
\shortauthorlist{Fernando A. Morales} 

\author{{
\sc Fernando A. Morales},\\[2pt]
Escuela de Matem\'aticas\\
Universidad Nacional de Colombia, Sede Medell\'in\\
Calle 50 A No 63--20, Bl 43, Of 106. Medell\'in, Colombia\\
{\email{famoralesj@unal.edu.co}} }

\maketitle

\begin{abstract}
{For a fissured medium, we analyze the impact that the geometry of the cracks, has in the phenomenon of preferential fluid flow. Using finite volume meshes we analyze the mechanical energy dissipation due to gravity, curvature of the surface and friction against its walls. We construct parameters depending on the \emph{Geometry} of the surface which are not valid for direct quantitative purposes, but are reliable for relative comparison of mechanical energy dissipation. Such analysis yields \emph{Information} about the preferential flow directions of the medium which, in most of the cases is not deterministic, therefore the respective \emph{Probability Spaces} are introduced. Finally, we present the concept of \emph{Entropy} linked to the geometry of the surface. This notion follows naturally from the random nature of the \emph{Preferential Flow Information}. }
{fissured media, energy dissipation, preferential flow, probability measures, geometric entropy.}
\\
2000 Math Subject Classification: 76S05, 97M99, 94A17
\end{abstract}

\section{Introduction}\label{Sec introl}
It is observed from experience, that the phenomenon of fluid flow through porous media is not uniform in every direction, on the contrary, preferential paths are developed. The problem of preferential flow has been extensively studied in recent years from several points of view and at different scales of modeling, due to its remarkable importance in different fields such as oil extraction, water supply, pollution of subsurface streams and soils, waste management, etc. At the pore scale, the presence of solutes and colloids, chemical reactions, high viscosity of the fluid and saturation level have been included in different theoretical and/or empirical models; see \cite{JensenHansenMagid, HagerdornMohn}. Nevertheless, upscaling this effects to that of the geological medium, or \emph{field scale}, has proved to be an extremely difficult task. A different approach emphasizes on the multiple scale aspects of the problem, when preferential flow occurs because of the presence of large pores (connected or not) or geological strata; which generates regions of fast and slow flow exchanging fluid. Hence, several models of coupled systems of partial differential equations have been proposed, such as dual \cite{VisaShow, Barenblatt, ArbogastDouglasHornung} and multiple porosity models \cite{SpagnuoloWright}, microstructure models \cite{Show97} and the coupling of laws at different scale: see \cite{ArbLehr2006} for an analytic approach of a Darcy-Stokes system and see \cite{LesinigoAngeloQuarteroni} for a numerical treatment to a Darcy-Brinkman system. On a different line there are several works, numerical \cite{Mclaren, ArbBrunson2007} or numerical/analytical \cite{WangFengHan}, dealing with the discretization and numerical aspects as well as the assessment and simulation of the proposed models. Yet another approach lies on a probabilistic point of view \cite{FioriJankovic, BerylandGolden, Golden} based on principles of conductivity. In this work we focus on the impact that the \emph{Geometry} of the cracks in a fissured medium has on this phenomenon. 
\begin{figure}
        \centering
        \begin{subfigure}[Fissured Medium]
                {\resizebox{6.5cm}{6.5cm}
{\includegraphics{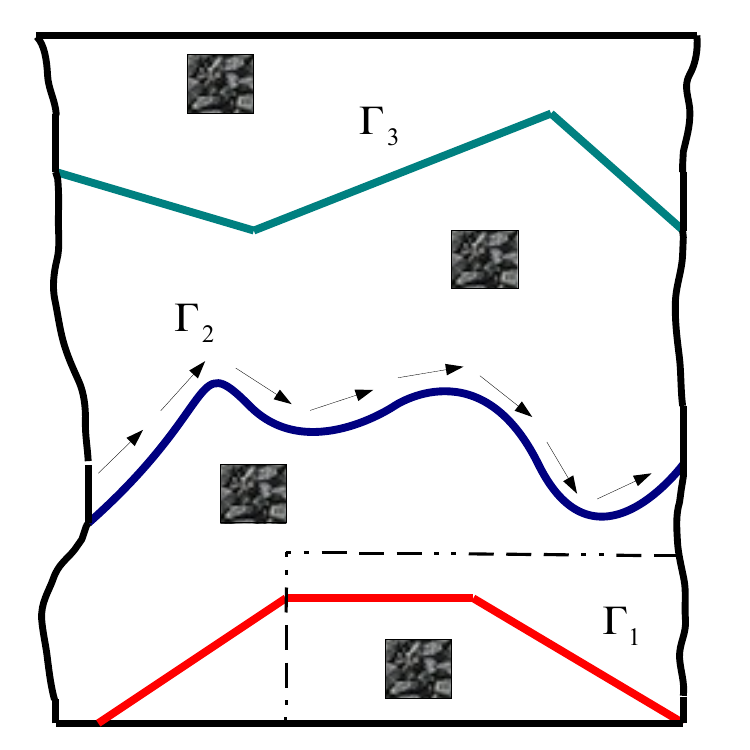} } }
                \label{Fig Geological Fissured Medium}
        \end{subfigure}
        \qquad
        ~ 
          \begin{subfigure}[Average Velocity on Fissured System Region]
                {\resizebox{6.5cm}{6.5cm}
{\includegraphics{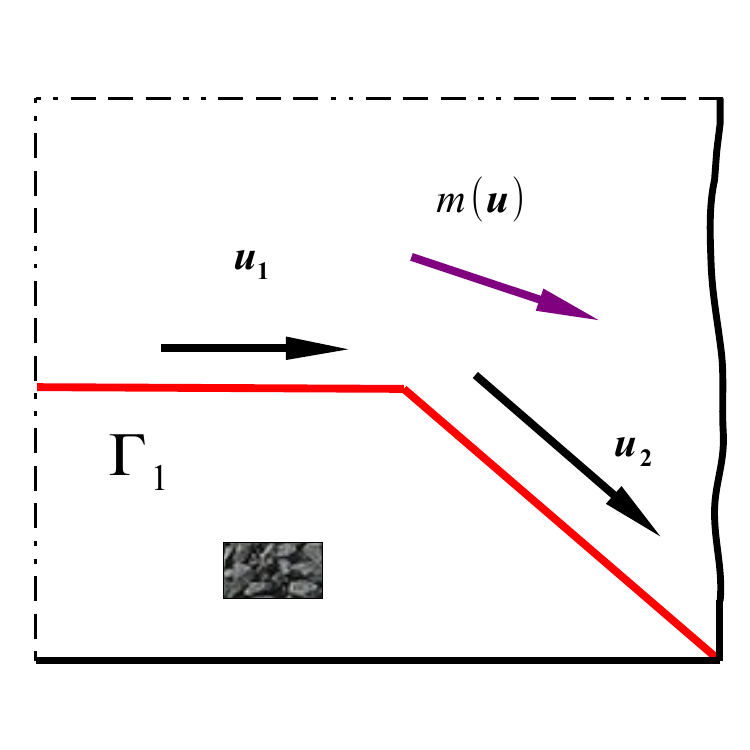} } }
\label{Fig Average Velocity on Fissure Portion}
        \end{subfigure}%
        ~ 
\caption{Fissured Medium and Average Velocities}
\end{figure}\label{Fig Fissured System}

Fissured media are common geological structures, here the fast flow occurs on the cracks while the rock matrix constitutes the slow flow region. For small values of the Reynolds number it is intuitive to see that the saturated flow on the fissures is predominantly parallel to the surface hosting it, see figure \ref{Fig Fissured System} (a). This fact has been shown in several rigorous mathematical works, see \cite{ArbogastDouglasHornung, WalkingtonShowalter} for homogenization techniques and \cite{SpagnuoloCoffield, Morales, MoralesShow2} for asymptotic analysis.  In the present work we exploit this fact to compute the ``average direction'' of tangential velocity fields hosted on a surface to predict the \emph{likely preferential flow directions} of the system. In figure \ref{Fig Fissured System} (b) a region of the fissured system is isolated in a way that contains only one crack. Assuming the flow is isotropic on the rock matrix, it follows that the preferential flow direction in this region strictly depends on the flow field hosted on the surface, i.e. its ``average'' behavior. In this work we assume the medium has one single fissure $\Gamma$ and state that its preferential flow is the ``average" tangential flow hosted on the surface $\Gamma$ of the crack. Towards the end of the exposition (section \ref{Closing Remarks}), it will be clear how to extend the method to a system with multiple fissures. 

Describing the exact flow field on the fissures on one hand, has a complexity level essentially equivalent to solving the problem of preferential direction itself, on the other hand for computational purposes it is always necessary to discretize the surface together with the flow fields. Consequently we choose a different approach: first we assume the medium has one single car we construct idealized flow fields related to finite triangulations of a surface, defined only by certain aspects of the geometry of the surface. Such fields are not realistic for describing the flow configuration on the manifold, however they are useful to compare quantities of mechanical energy dissipation. Hence, in order to find the preferential flow directions amongst all the possible aforementioned flow fields, we apply the reasoning line of the Arquimedian weight comparison method. As it turns out, in most of the cases the preferential flow direction is not unique, and a probabilist treatment must be adopted.

Next, we introduce the notation, vectors in $\R^{\!2}$ and $\R^{\!3}$ will be denoted with bold characters and $\vert \cdot \vert$ stands for the Euclidean norm. If $\mathbf{x} = (x_{1}, x_{2}, x_{3})\in \R^{\!3}$ we denote $\widetilde{x} = (x_{1}, x_{2})$, and $\mathbf{x} = (\widetilde{x}, x_{3})$. The notation $m(\cdot)$ indicates computation of average according to the context. The unitary circle in $\R^{\!2}$ is denoted by $S^{1}$ and the unitary sphere in $\R^{!3}$ is denoted by $S^{2}$. For a given set $A$ we denote $\vert A \vert$ its 2-D or 1-D Lebesgue measure in both cases, since it will be clear from the context and $\#A$ stands for its cardinal. The paper is organized as follows: in sections 2, 3 and 4 the preferential flow phenomenon is analyzed from the perspectives of Curvature, Gravity and Friction respectively. All of them are based on comparing mechanical energy dissipation, therefore each section begins constructing adequate flow fields to quantify these losses. Next, the exposition moves to a rigorous discussion of mathematical aspects of the model; most of them necessary for a successful construction of the \emph{Preferential Flow Directions Probability Space}. Section 5 shows how to assemble the effects previously analyzed and defines the entropy of the preferential flow information; then closes with final remarks and future work. In the reminder of this section we present the geometric setting and the minimum necessary background from fluid mechanics.   
%
%
\subsection{Geometric Setting and Triangulation of the Surface}\label{Sec geometric setting and triangulation}
\begin{figure}[!]
\caption[16]{Edge Influence $\&$ Configuration}\label{Fig Edge Influence}
\includegraphics{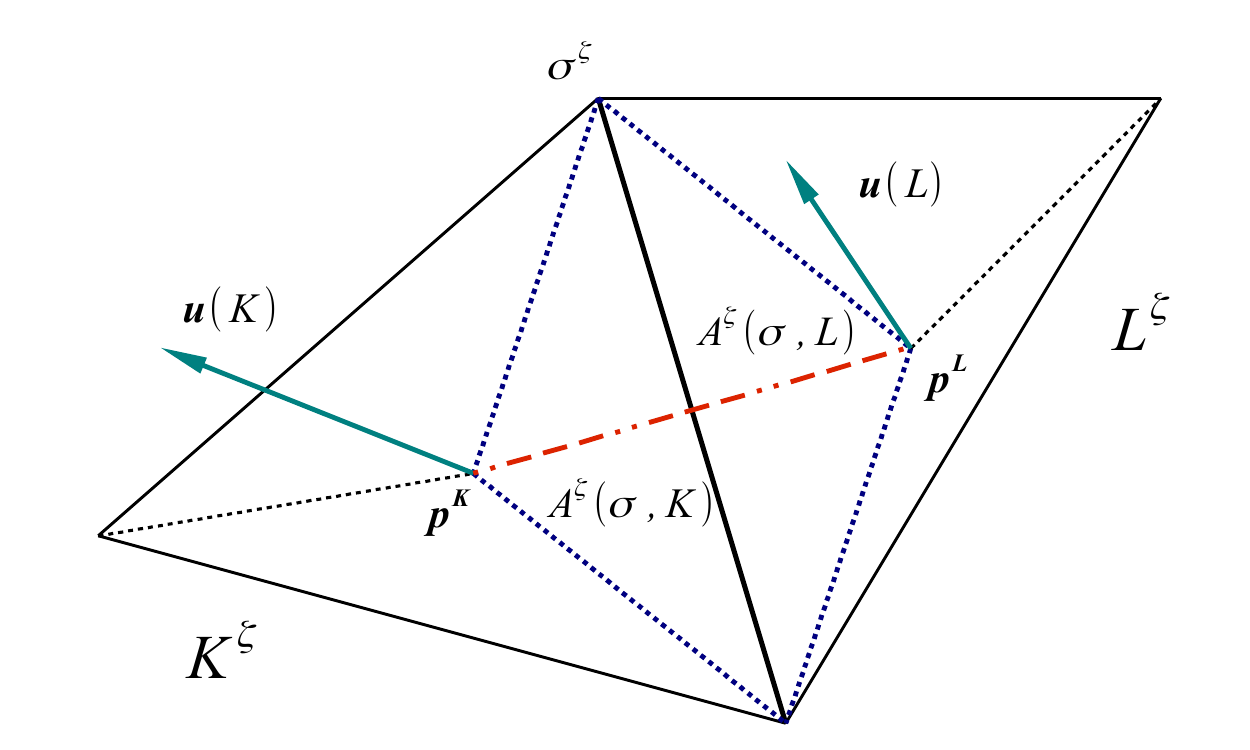}
\end{figure}

This work will be restricted to the treatment of surfaces $\Gamma$ coming from a piecewise $C^{1}$-function $\zeta:G\rightarrow \R$ defined on an open bounded simply connected set $G$ of $\R^{\!2}$. In particular the piecewise $C^{1}$ 2-D manifold $\Gamma$ has an Atlas containing one element. The approximation of surfaces will be done by piecewise linear affine triangulations. However, the triangulation has to meet certain geometric conditions in order to be suitable for later quantifications of mechanical energy losses. Those conditions are consistent with the concept of \emph{admissible mesh} in the sense of \cite{GallouetHerbin, BradjiHerbin} (Definition 9.1 page 762); we have
\begin{definition}\label{Def admissible mesh}
Let $\Omega$ be an open bounded polygonal subset of $\R^{d}$, $d = 2, 3$. An admissible finite volume mesh of $\Omega$, denoted by $\triang$, is given by a family of ``control volumes" $\{K: K\in \triang\}$, which are open polygonal convex subsets of $\Omega$, a family of subsets of $cl(\Omega)$ contained in hyperplanes of $\R^{d}$, denoted by $\E$ (these are the edges in two dimensions or faces in three dimensions of the control volumes), with strictly positive $(d-1)$-dimensional measure, and a family of points, of $\Omega$ denoted by $\control$ satisfying the following properties:
\begin{enumerate}[(i)]
\item
$ cl(\Omega) = cl \bigcup_{K\in\mathcal{T}} K$

\item
For any $K\in \triang$, there exists a subset $\E_{K}$ of $\E$ such that $\partial K = cl(K)-K = \bigcup_{\sigma\in \E_{K}} cl(\sigma)$. Furthermore, $\E = \bigcup_{K\in \mathcal{T}} \E_{K}$.

\item 
For any $(K, L)\in\triang^{\,2}$ with $K\neq L$, either the $(d-1)$-dimensional Lebesgue measure of $cl(K)\cap cl(L)$ is 0 or $cl(K)\cap cl(L) =cl(\sigma)$ for some $\sigma\in \E$, which will then be denoted by $K\vert L$.

\item
The family $\control = \{x_{K}: K\in \triang\}$ is such that $x_{K}\in cl(K)$ (for all $K\in \triang$) and if $\sigma = K\vert L$, it is assumed that $x_{K}\neq x_{L}$, and that the straight line $\D_{K, L}$ going through $x_{K}$ and $x_{L}$ is orthogonal to $K\vert L$.

\item
For any $\sigma\in \E$ such that $\sigma\subset\partial \Omega$, let $K$ be the control volume such that $\sigma\in \E_{K}$. If $x_{K}\notin \sigma$, let $\D_{K, \sigma}$ be the straight line going through $x_{K}$ and orthogonal to $\sigma$, then the condition $\D_{K, \sigma}\cap \sigma\neq \emptyset$ is assumed. Define 
\begin{equation}\label{Def transimition point on 2-D mesh}
   y_{\sigma} \defining  \D_{ K, \sigma}\cap \sigma
\end{equation}
\end{enumerate}
\end{definition}
From now on we adopt triangular meshes as provided in \cite{GallouetHerbin}
\begin{definition}\label{Def admissible triangular mesh}
Let $\Omega$ be an open bounded polygonal subset of $\R^{2}$. We say a triangular mesh is a family $\triang$ of open triangular disjoint subsets of $\Omega$ such that two triangles having a common edge have also two common vertices and such that all the interior angles of the triangles are less than $\frac{\pi}{2}$. 
\end{definition}
Clearly a triangular mesh described in the definition above meets the conditions of \ref{Def admissible mesh}. In particular, the condition on the interior angles assures that the orthogonal bisectors intersect inside each triangle, thus naturally defining the points $x_{K}\in K$. Since definitions \ref{Def admissible mesh} and \ref{Def admissible triangular mesh} demand a polygonal domain we introduce the collection of eligible polygons.
\begin{definition}\label{Def eligible polygons}
Let $\Gamma = \{[\widetilde{x}, \zeta(\widetilde{x})]: \widetilde{x}\in G\}$ be a piecewise $C^{1}$ surface. We say a polygonal domain is eligible for triangulation of $\Gamma$ if it is contained in $G$ and if its vertices lie on the boundary of $G$. From now on we denote $Poly\, (G)$ the family of all such polygons.
\end{definition}
Now we introduce a central definition for the type of triangulations to be worked on
\begin{definition}\label{Def lifting}
Let $\Gamma = \{[\widetilde{x}, \zeta(\widetilde{x})]: \widetilde{x}\in G\}$ be a piecewise $C^{1}$ surface and $K$ a triangular domain contained in $G$ with vertices $\{z_{\,\ell}: 1\leq \ell\leq 3\}\subset \R^{2}$ we define its ``lifting" as the closed convex hull of the points $\{[z_{\,\ell}, \zeta(z_{\,\ell})]:1\leq \ell\leq 3\}\subset \R^{3}$. We denote this surface by $\Kz$ and the outer unitary vector perpendicular to it by $\n(K)$.
\end{definition}
Next we define a \emph{Triangulation} of the surface $\Gamma$.
\begin{definition}\label{Def triangulation of gamma}
Let $\Gamma$ be a piecewise $C^{1}$ surface, $\Hp\in Poly\, (G)$ and $\triang$ an admissible triangular mesh of $\Hp$ as in definition \ref{Def admissible triangular mesh}. 
\begin{enumerate}[(i)]
   \item 
   We say the triangulation of $\Gamma$ relative to the polygon $\Hp$ and the mesh $\triang$, is given by the ``lifting" $\Kz$ of each element $K$ of$\triang$. We denote
\begin{equation}\label{Def triangulation surface P, T}
\Gamma_{\scriptscriptstyle \Hp, \,\triang}\defining\bigcup\left\{\Kz: K\in \triang\right\}
\end{equation}

\item
The point of control $\pk$ is given by the unique point in $\Kz$ such that its horizontal projection agrees with $x_{K}$, the circumcenter of $K$ i.e.
\begin{equation}\label{Def triangulation points of control P, T}
\pk - ( \pk\cdot\kversor ) \, \kversor = x_{K}
\end{equation}
Moreover $\pk$ is the ``lifting'' of $x_{K}$.  

\item
Given an edge $\sigma \in \E$ and its middle point $y_{\sigma} $, the associated ``transmission point'' is the unique point $\qs$ contained in $\Gamma_{\Hp, \triang}$ such that
\begin{equation}\label{Def triangulation points of transmission}
\qs -( \qs\cdot\kversor ) \, \kversor = y_{\sigma}
\end{equation}
i.e. $\qs$ is the ``lifting'' of $y_{\sigma}$. 

\item
Define $\E_{int}\defining\{\sigma\in \E: \vert \sigma\cap \partial\Hp\vert = 0 \}$ i.e. the set of interior edges of the triangulation $\triang$.

\item
For each element $K\in \triang$ define its \emph{edge-influence triangles} as the three subtriangles generated by drawing rays from $x_{K}$ to each of its vertices. Figure \ref{Fig Edge Influence} depicts the lifting of two neighboring elements $K$ and $L$, its common edge $\sigma = K | L$ and the corresponding edge-influence triangles.

\end{enumerate}
\end{definition}
%
   %
   %
%
%
%
%
%
%
%
\subsection{The Strain Rate Tensor and Mechanical Energy Loss}\label{Sec strain rate tensor}
For the sake of completeness we recall the definition of strain rate tensor \cite{Batchelor}. Given a differentiable flow field $\v: G\rightarrow \R^{d}$, $G$ open set in $\R^{d}$, the strain rate tensor $D(\v): G\rightarrow \R^{d\times d}$ is given by
\begin{equation}\label{Eq Strain Rate Tensor}
D_{j, \ell}(\v) \defining \frac{1}{2}\left(\frac{\partial \v_{j}}{\partial x_{\ell}}
+ \frac{\partial \v_{\ell}}{\partial x_{j}}\right), \quad 1\leq j, \, \ell\leq d
\end{equation}
Finally, the internal deformation energy of a viscous fluid is given by \cite{Bear}
\begin{equation}\label{Eq dissipation of mechanical energy}
E = \frac{2\,\mu}{\rho}\,D(\v) : D(\v) \defining 
\frac{2\,\mu}{\rho}\sum_{j, \ell}\left\vert D_{j, \ell}(\v)\right \vert^{2}
\end{equation}
Where $\mu$ represents the viscosity and $\rho$ the density of the fluid. 
%
%
%
%
\section{Preferential Flow due to Curvature}\label{Sec curvature dissipation C2 general}
%
%
\subsection{Flow Hypothesis}\label{Sec flow hypothesis}
We want to compute a conservative tangential flow field, hosted within the surface $\Gamma$ and totally defined by its curvature. As already specified in the introduction, this paper will be restricted to the construction of a discretized flow field related to a triangulation $\Gamma_{\scriptscriptstyle \Hp, \,\triang}$. The changes on the flow field must be exclusively  due to the changes of directions on the elements of the surface $\Gamma_{\scriptscriptstyle \Hp, \,\triang}$. Then, for simplicity we choose the following defining properties
\begin{enumerate}[(i)]\label{Def conditions on the artificial flow field}

\item
The velocity must be constant in magnitude and direction within a flat face.

\item The magnitude of the velocity must be constant on every part of the surface.

\item
The field must meet the continuity flow condition i.e. on the edge where two different faces intersect the component of the velocities perpendicular to the edge must have the same magnitude.
\end{enumerate}
For the construction of such velocity field we introduce a velocity of reference or \emph{master velocity} which will be denoted $\ugamma$. Since the surface $\Gamma$ is defined by a $C^{1}(G)$ function, no triangulation $\Gamma_{\Hp,\,\triang}$ contains vertical faces, i.e. $\n (K)\cdot\kversor\neq 0$ for all $K\in \triang$. Hence, whichever tangential flow that the surface $\Gamma_{\Hp,\,\triang}$ hosts has a non-null projection onto the plane $\langle \kversor\rangle^{\perp}$. Consequently, it is enough to assume that the velocity of reference $\ugamma$ is hosted in the horizontal plane. 
%
%
%
%
\subsection{Construction of the Velocity Field}\label{Sec velocity field construction}
Let $\Gamma$ be a piecewise $C^{1}$ surface and $\Gamma_{\scriptscriptstyle \Hp, \triang}$ be a triangulation. Given a reference velocity $\ugamma$ we are to build the velocity $\u(K)$ on the element $\Kz$(the lifting of $K$). If $\Kz$ is horizontal i.e if $\n(K) \equiv \kversor$ we simply set $\u(K) = \ugamma$. For the non-trivial case when $\Kz$ is not horizontal ($\n(K)\times \kversor\neq 0$), we proceed as follows. In the figure \ref{Fig Velocities Relation} below we illustrate the relation between velocities. It depicts the horizontal and vertical view of the intersection between the plane $\langle \kversor\rangle^{\perp}$ and the plane containing $\Kz$ an element of $\Gamma_{\Hp,\,\triang}$. 
%
%
\begin{figure}
        \centering
        \begin{subfigure}[Horizontal View]
                {\resizebox{6.5cm}{6.5cm}
{\includegraphics{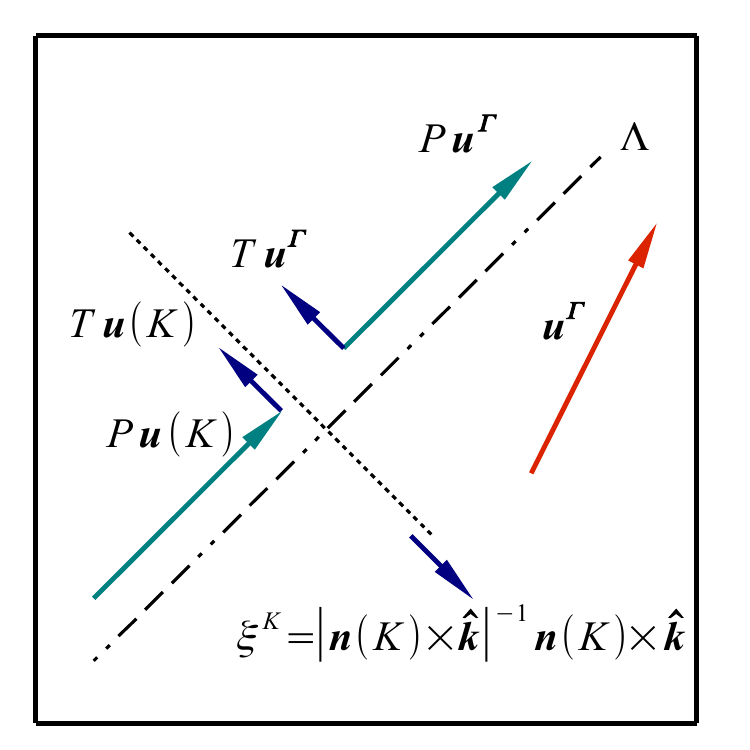} } }
                \label{Fig Velocities Relation Horizontal View}
        \end{subfigure}
        \qquad
        ~ 
          \begin{subfigure}[Vertical View]
                {\resizebox{6.5cm}{6.5cm}
{\includegraphics{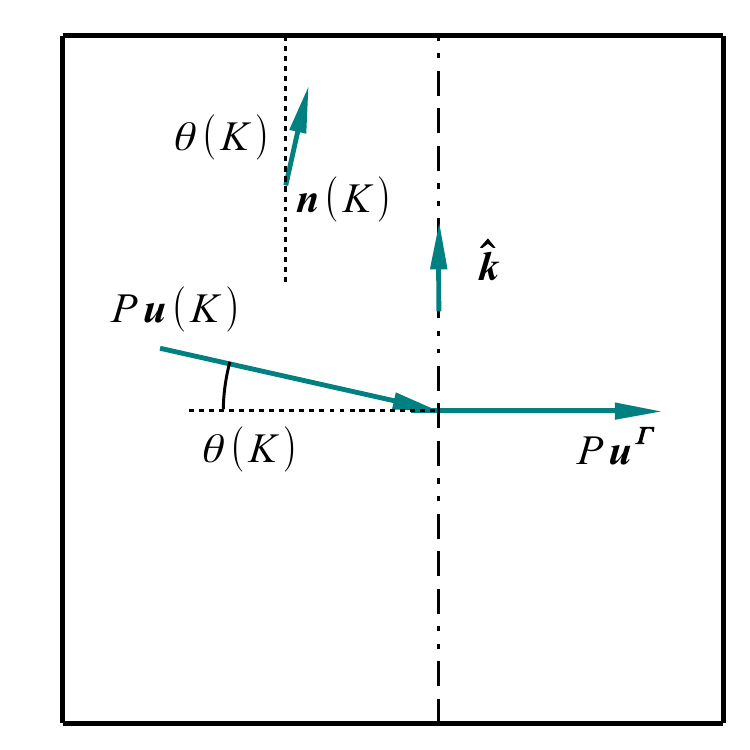} } }
\label{Fig Velocities Relation Vertical View}
        \end{subfigure}%
        ~ 
\caption{Intersection $\Kz\cap \langle \kversor \rangle^{\perp}$}
\end{figure}\label{Fig Velocities Relation}
On the left hand side of figure \ref{Fig Velocities Relation} we have the reference velocity $\ugamma$  and a decomposition of the local velocity $\u(K)$. The fine dashed line in the direction of the unitary vector $\xi^{K}\defining\vert\n(K)\times \kversor\vert^{-1}\,\n(K)\times \kversor$ represents the intersection line of the plane containing $\Kz$ and the plane $\langle\kversor\rangle^{\perp}$. We decompose $\ugamma$ in two vectors lying on the horizontal plane $\langle\kversor\rangle^{\perp}$, one  parallel to the intersection line $\langle \xi^{K}\rangle$ and the other perpendicular to it i.e.
\begin{subequations}\label{Def velocity projections}
\begin{equation}\label{Def velocity parallel to the intersection}
T(K)\, \ugamma = T\, \ugamma\defining \left(\ugamma\cdot\xi^{\,K}\right)\xi^{\,K} ,\\
%
\end{equation}
\begin{equation}\label{Def velocity perpendicular to intersection}
P(K)\,\ugamma = P\, \ugamma\defining
\ugamma - \left(\ugamma\cdot\xi^{\,K}\right)\xi^{\,K} ,
\end{equation}
\begin{equation}\label{Def velocity decomposition}
\ugamma = T\, \ugamma +
P\,\ugamma .
\end{equation}
\end{subequations}
Since the component $T\,\ugamma$ belongs to the intersection of both planes $\langle\kversor\rangle^{\perp}\cap \langle\n (K)\rangle^{\perp}$, we set it equal to the component of $\u(K)$ in the same direction i.e.
%
%
%
\begin{equation}\label{Def velocity intersection line projection}
T(K)\, \u(K) \defining T(K)\, \ugamma
%
\end{equation}
%
%
%
On the right hand side of figure \ref{Fig Velocities Relation} we depict the trace through the vertical plane $\Lambda \defining \langle\xi^{\,K}\rangle^{\perp}$. Denote $P(K)\, \u(K)$ the component of $\u(K)$ perpendicular to $\xi^{\,K}$, we set this component to be a rotation of $P(K)\,\ugamma$ by the angle $\theta(K)$, strictly contained in the plane $\langle\xi^{\,K}\rangle^{\perp}$; where the angle $\theta(K)$ is equal to the angle formed between $\kversor$ and $\n (K)$.
%
%
%
%
Notice that the map $\ugamma \in \R^{\!2}\mapsto \u(K) $ for fixed $\Kz\in\Gamma_{\scriptscriptstyle \Hp, \,\triang}$ is linear, therefore writing $\ugamma = \alpha_{\,1}\,\iversor+ \alpha_{\,2}\,\jversor$ a direct calculation yields
\begin{subequations}\label{Eq velocity due curvature 3-D}
\begin{equation}\label{Eq local velocity due curvature 3-D}
 \u(K)
\equiv \left(\begin{array}{cc}
1- \dfrac{\n_{1}^{\,2}}{1+\n_{\,3}}   & -\dfrac{\n_{1}\,\n_{\,2}}{1+\n_{\,3}}\\[8pt]
-\dfrac{\n_{1}\,\n_{\,2}}{1+\n_{\,3}} & 1- \dfrac{\n_{\,2}^{\,2}}{1+\n_{\,3}}\\[8pt]
\n_{1} & \n_{\,2}
\end{array}
\right)
\left\{\begin{array}{c}
\alpha_{\,1}\\[20pt]
\alpha_{\,2}
\end{array}
\right\}\\
\defining 
\left(V_{\Gamma_{\Hp,\,\triang}}\ugamma \right)(K) .
\end{equation}
Where $\n(K) = \n_{1} \iversor +  \n_{2}\jversor + \n_{3} \kversor$ . The global velocity field $V_{ \Gamma_{ \Hp,\,\triang}}\ugamma : \Hp\rightarrow \R^{\!3}$ is defined by
\begin{equation}\label{Eq global velocity due curvature 3-D}
 V_{ \Gamma_{ \Hp,\,\triang}}\ugamma(s) \defining
 \sum_{K\in \triang} \u(K) \, \ind_{K}(s)\,,\quad s\in \Hp
\end{equation}
\end{subequations}

\begin{remark}\label{Rem comments on the velocity field}
\begin{enumerate}[(i)]
\item
By construction it is clear that the field in the expression \eqref{Eq velocity due curvature 3-D} meets the criteria (i), (ii) and (iii) of the flow hypothesis described at the begining of section \ref{Sec flow hypothesis}.

\item
Observe that if $\Kz$ is horizontal, then $\u(K) = \ugamma$ as expected i.e. the expression \eqref{Eq velocity due curvature 3-D} covers all the possible cases. 

\item
Whenever there is no ambiguity we will simply denote $\u = V_{ \Gamma_{ \Hp,\,\triang}}\ugamma$.
\end{enumerate}
\end{remark}
For the field of velocity given by expression \eqref{Eq velocity due curvature 3-D} we define the average velocity in the natural way
\begin{equation}\label{Def average velocity operator}
m( V_{\Gamma_{\Hp,\,\triang}}\b  )\defining 
%
\sum_{K\,\in\,\triang} \frac{\vert\Kz \vert}{\vert\Gamma_{\Hp,\,\triang}\vert}
\, \left( V_{\Gamma_{\Hp,\,\triang}}\b \right) (K)
\end{equation}
Clearly, the map $m\circ  V_{\Gamma_{\Hp,\,\triang}}: \R^{2}\rightarrow \R^{3}$ is linear. It is important to stress that  $m( V_{\Gamma_{\Hp,\,\triang}}\b )$ may not be tangential to (or hosted within) $\Gamma_{\Hp,\,\triang}$. Another important fact is the following 
\begin{lemma}\label{Th average velocity operator}
   Let $\Gamma$ be a piecewise $C^{1}$ surface, $\Gamma_{\scriptscriptstyle \Hp, \triang}$ be a triangulation and the average velocity operator $m$ defined by \eqref{Def average velocity operator}, then
   \begin{enumerate}[(i)]
   \item \label{Th null kernell average velocity operator}
   $\ker(m\circ  V_{\Gamma_{\Hp,\,\triang}}) = \{\zero\}$. 
   
   \item \label{Th rank average velocity operator}
   The space $(m\circ  V_{\Gamma_{\Hp,\,\triang}})(\R^{\!2})$ is two dimensional.  
   \end{enumerate}
   \begin{proof}
   \begin{enumerate}[(i)]
      \item
      Let $\b = \alpha_{1}\iversor + \alpha_{2} \jversor$ such that $\vert \b \vert = 1$ and $K\in \triang$ be arbitrary, then
      \begin{multline*}
         \left(V_{\Gamma_{\Hp,\,\triang}}\b\right)(K)\cdot \b  = \left(1 - \frac{\n_{1}^{\,2}}{1+\n_{\,3}}\right) \alpha_{1}^{2}
         - 2\,\frac{\n_{1}\, \n_{2}}{1+\n_{\,3}}\,\alpha_{1}\, \alpha_{2}
         + \left(1 - \frac{\n_{2}^{\,2}}{1+\n_{\,3}} \right)\alpha_{2}^{2} \\
         = 
         (\alpha_{1}^{2} + \alpha_{2}^{2}) - \frac{1}{1 + \n_{\,3}} (\n_{1} \alpha_{1} + \n_{2}\alpha_{2})^{2} \geq 1 - (\vert \b\vert \vert \n\vert)^{2} = 0
      \end{multline*}
      Therefore, the horizontal projection of $\left(V_{\Gamma_{\Hp,\,\triang}}\b\right)(K)$ makes an angle with $\b$ less or equal than $\frac{\pi}{2}$. This implies that the projection onto $\b$ satisfies
      \begin{equation*}
         P_{\scriptscriptstyle \b} \, \left(V_{\Gamma_{\Hp,\,\triang}}\b\right)(K) = \lambda \, \b \, , \quad \lambda > 0\, ,\; \forall\, K\in \triang .
      \end{equation*}
      Finally, since $\b\neq 0$ and the weighting coefficients in \eqref{Def average velocity operator} multiplying $\left(V_{\Gamma_{\Hp,\,\triang}}\b\right)(K)$ are positive for all $K\in \triang$, it follows that $P_{\scriptscriptstyle \b} \left(m_{\scriptscriptstyle \Gamma_{\Hp,\,\triang}}(\b)\right) \neq \zero$; which concludes the first part.
      
      \item
      Follows immediately from the previous part and the dimension theorem.
      \end{enumerate}
   \end{proof}
\end{lemma}
%
%
%
%
%
%
%
%
%
\subsection{Dissipation of Mechanical Energy Model due to Curvature}\label{Sec dissipation due to curvature}
%
%
%
%
The discrete model of mechanical energy dissipation due to change of direction has to be consistent with the expression \eqref{Eq dissipation of mechanical energy} i.e. we need to generate a discrete field of \emph{strain rate tensors} using the velocity given in \eqref{Eq velocity due curvature 3-D}. The flow field can change only from one element of the triangulation to another. Consequently, the variations of the flow field across the edges define the strain rate tensor we seek. 
\begin{definition}
Let $\sigma\in \E_{int}$ and $K$, $L$ be the two elements of $\Gamma_{\scriptscriptstyle \Hp, \,\triang}$ such that $\sigma = K\vert L$; denote $\pk, \pl$ and $\u(K), \u(L)$ the respective points of control and fluid velocity  for each element. 
\begin{enumerate}[(i)]
\item
Define the strain rate tensor across $\sigma$ by

\begin{subequations}\label{Def discreet strain rate tensor}
\begin{equation}\label{Def discreet strain rate tensor 1}
D_{k,\,\ell}(\sigma)\defining \frac{1}{2}\,\frac{\left[\u(K)- \u(L)\right]\cdot\eversor_{k}}
{\left[\pk - \pl\right]\cdot\eversor_{\ell}}
+\frac{1}{2}\, \frac{\left[\u(K)- \u(L)\right]\cdot\eversor_{\ell}}
{\left[\pk - \pl\right]\cdot\eversor_{k}} \, ,
\end{equation}
\begin{equation}\label{Def discreet strain rate tensor 2}
D(\sigma)\defining \left\{D_{k,\,\ell}(\sigma):1\leq k, \ell\leq 3\right\} .
\end{equation}
\end{subequations}
Here $\eversor_{1} = \iversor$, $\eversor_{\,2}=\jversor$ and $\eversor_{3}= \kversor$.

\item
The strain tensor is extended to the whole surface $\Gamma_{\scriptscriptstyle \Hp, \,\triang}$ by
\begin{equation}\label{Def distribution of strain rate tensor}
D (s)\defining \sum_{\sigma\,\in\,\E_{int}}D (\sigma)
\left[\ind_{A(\sigma, K)}+\ind_{A(\sigma, L)}\right](s)\,,\quad s\in \Hp .
\end{equation}
Where $A(\sigma, K)$ and $A(\sigma, L)$ denote the edge-influence triangles of $K$ and $L$ respectively, incident on $\sigma$ given in (v) definition \ref{Def triangulation of gamma}. Figure \ref{Fig Edge Influence} displays a horizontal view of two neighboring elements. 
\end{enumerate}
\end{definition}
%
%
%
%
%
%
%
Finally, using \eqref{Def distribution of strain rate tensor} to compute \eqref{Eq dissipation of mechanical energy} we have that the global dissipation of energy on the triangulation $\Gamma_{\scriptscriptstyle \Hp, \,\triang}$ under the \emph{master velocity} $\ugamma$ is given by
\begin{equation}\label{Eq dissipation due to change of direction}
U_{curv} (\ugamma)\defining \frac{2\,\mu}{\rho}
\sum_{\sigma\,\in\,\E_{int}}
\left\{\vert A^{\zeta}(\sigma, K)\vert+\vert A^{\zeta} (\sigma, L)\vert\right\}
D (\sigma) : D (\sigma) .
\end{equation}
Here $\vert A^{\zeta}(\sigma, K)\vert, \vert A^{\zeta}(\sigma, L)\vert$ are the areas of the lifted adjacent triangles $A(\sigma, K)$, $A(\sigma, L)$ respectively. Clearly $U_{curv}$ depends only on the triangulation $\Gamma_{\Hp,\,\triang}$ and the master velocity $\ugamma$.
%
%
%
%
\begin{remark}\label{Rem Comments on the Mechanical Energy Dissipation due to Curvature}
Given a triangulation $\Gamma_{\scriptscriptstyle \Hp,\,\triang}$ of a piecewise $C^{1}$ surface $\Gamma$, denote by $diam(\Gamma_{\scriptscriptstyle  \Hp,\,\triang})$ the maximum diameter of its elements $\Kz$. Letting $diam(\Gamma_{\Hp,\,\triang})\rightarrow 0$, on one hand, the map $s\in \Hp \mapsto \sum_{K\in \triang} \n(K)\, \ind_{K}$ converges (non-conformally) to $s\in \Hp \mapsto \n(s)$ almost everywhere; on the other hand, the distance $\vert \pk - \pl \vert$ tends to zero. Due to the definition of the flow field $V_{\Gamma_{\Hp,\,\triang}}\ugamma$ given in \eqref{Eq velocity due curvature 3-D}, the expression \eqref{Def discreet strain rate tensor} starts approaching values of a directional derivative of the map $s\in \Gamma\mapsto \n(s)$ (on the points where is differentiable); and the curvature information of the surface $\Gamma$ is contained in these derivatives. Hence, the tensor proposed in \eqref{Def discreet strain rate tensor} and the mechanical energy dissipation proposed in \eqref{Eq dissipation due to change of direction} are heavily defined by the curvature (or rather an approximation of the curvature) of the surface $\Gamma$.
\end{remark}
%
%
%
%
%
\subsection{Minimum and Maximum Mechanical Energy Dissipation due to Curvature}\label{Sec mininimal dissipation due to curvature}
By definition the functional $U_{curv}: \R^{2}\rightarrow \R$ is a quadratic form, then it holds that
\begin{equation}\label{Eq dissipation due to change of direction matrix}
U_{curv} (\b) = \b^{T}\,M_{curv}\,\b ,
\end{equation}
for $M_{curv}$ symmetric, positive semi-definite matrix. Due to the spectral theorem, the matrix $M_{curv}$ is orthogonally diagonalizable. Denote $0\leq \lambda_{1}\leq \lambda_{\,2}$ the eigenvalues and $\{\fversor_{\!\!\!1}, \fversor_{\!\!2}\}$ an associated orthonormal basis of eigenvectors, then 
\begin{equation}
\begin{split}
U_{curv}(\b) & = \lambda_{1}\, (\b\cdot\fversor_{\!\!\! 1})^{2} + 
\lambda_{\,2} (\b\cdot\fversor_{\!\!2} )^{2} , \\
\lambda_{1} & = \min\{ U_{curv}(\b):\vert\b\vert = 1\} = U_{curv}(\fversor_{\!\!\! 1}) , \\
\lambda_{2} & = \max\{ U_{curv}(\b):\vert\b\vert = 1\} = U_{curv}(\fversor_{\!\!\! 2}) .
\end{split}
\end{equation}
 i.e. the question of minimum and maximum mechanical energy dissipation due to curvature of the surface is equivalent to an eigenvalue problem of $M_{curv}\in \R^{2\times 2}$.
%
%
%
%
\subsection{Preferential Fluid Flow Directions Due to Curvature and Their Probability Space.}\label{Sec preferential direction due to curvature}
The \emph{Preferential Fluid Flow Directions} of the surface $\Gamma_{\Hp,\,\triang}$ due to \emph{Curvature} are given by
\begin{equation}\label{Def curvature preferential direction}
\varpi_{curv} \defining \left\{\frac{m(V_{\Gamma_{\Hp,\,\triang} }\b)}{\vert m(V_{\Gamma_{\Hp,\,\triang} }\b)  \vert}: \b\neq \zero, \,U_{curv}(\b) = \lambda_{1}\vert\b\vert ^{\,2}\right\} .
\end{equation}
%
%
Due to lemma \ref{Th average velocity operator} part \ref{Th null kernell average velocity operator} the set $\varpi_{curv}$ is well-defined. It is direct to see that if $\lambda_{1}<\lambda_{2}$ then $\varpi_{curv}$ will have two elements, namely $\vert m_{\scriptscriptstyle \Gamma_{\Hp,\,\triang}}(\fversor_{\!\!\! 1})\vert^{-1} m_{\scriptscriptstyle \Gamma_{\Hp,\,\triang}}(\fversor_{\!\!\! 1})$ and $\vert m_{\Gamma_{\Hp,\,\triang}}(-\fversor_{\!\!\! 1})\vert^{-1} m_{\Gamma_{\Hp,\,\triang}}(-\fversor_{\!\!\! 1})$ for $\fversor_{\!\!\! 1}$ the unitary vector associated to $\lambda_{1}$. On the other hand if $\lambda_{1} = \lambda_{2}$ then $\varpi_{curv} $ has infinitely many elements due to lemma \ref{Th average velocity operator} part \ref{Th rank average velocity operator}. In both cases it can not be chosen which direction within $\varpi_{curv} $ is preferential over the others. Due to the uncertainty of this information we must treat it from a \emph{Probabilistic} point of view. In order to give a consistent definition for the probability space of preferential directions we need to introduce a previous one for technical reasons.
\begin{definition}\label{Def sigma algebra curvature}
Let $\Gamma$ be a piecewise $C^{1}$ surface, $\Gamma_{\scriptscriptstyle \Hp, \triang}$ be a triangulation and $\omega_{curv} = \{\b\in S^{1}: U_{curv} (\b) = \lambda_{1}\}$, consider the surjective function
\begin{equation}\label{Def probabilistic distribution function curvature}
\varphi:\omega_{curv} \rightarrow \varpi_{curv}\,,\quad
\varphi(\eversor)\defining \frac{m(V_{\Gamma_{\Hp,\,\triang} }\eversor)}{\vert m(V_{\Gamma_{\Hp,\,\triang} }\eversor) \vert} .
\end{equation}
Let $\beta$ be the family of all Borel sets of $S^{1}$ intersected with $\omega_{curv}$, define the following $\sigma$-algebra
\begin{equation}\label{Eq sigma algebra curvature}
\tau_{curv}\defining
\left\{A \in \wp(\varpi_{curv}): \varphi^{-1}(A) \in \beta\right\}.
\end{equation}
Where $\wp(\varpi_{curv})$ is the power set of $\varpi_{curv}$.
\end{definition}
Finally, we endow the preferential flow space with the uniform probability distribution. 
\begin{definition}\label{Def distribuition due to curvature}
Let $\Gamma$ be a piecewise $C^{1}$ surface, $\Gamma_{\scriptscriptstyle \Hp, \triang}$ be a triangulation and $\varpi_{curv}$ be the associated preferential fluid flow directions defined in \eqref{Def curvature preferential direction} then
\begin{enumerate}[(i)]
\item
If $\lambda_{1}<\lambda_{\,2}$ then $\# \varpi_{curv} = 2$; define
\begin{equation}\label{Def dicrete probabilistic distribution curvature}
\prob_{curv}\left\{\frac{m(V_{\Gamma_{\Hp,\,\triang} }\fversor_{\!\!\! 1} ) }{\vert m(V_{\Gamma_{\Hp,\,\triang} }\fversor_{\!\!\!  1})\vert}\right\} \defining \frac{1}{2}\,,\quad
\prob_{curv}\left\{\frac{m( V_{\Gamma_{\Hp,\,\triang} }(-\fversor_{ \!\!\! 1} ) ) }{\vert m(V_{\Gamma_{\Hp,\,\triang} }(-\fversor_{ \!\!\! 1}) )\vert}\right\} \defining \frac{1}{2}.
\end{equation}
Where $\fversor_{\!\!\!1}$ is the eigenvector associated to $\lambda_{1}$. 

\item
If $\lambda_{1} = \lambda_{\,2}$ define
\begin{equation}\label{Def continuous probabilistic distribution curvature}
   \prob_{curv}[A] \defining \frac{1}{\mu (\varphi^{-1} (\varpi_{curv} ) ) }
   \, \mu (\varphi^{-1} ( A ) ) = \frac{1}{\mu (\omega_{curv} )  }
   \, \mu (\varphi^{-1} ( A ) )\, , \quad \forall \, A \in \tau_{curv}.
\end{equation}
Here $\mu$ indicates the arc-length measure in $S^{1}$.
\end{enumerate}
\end{definition}
%
%
%
%
%
%
%
%
%
%
%
%
%
\section{Preferential Flow Due to Gravity}\label{Sec preferential flow due to gravity}
In the present section we address the question of quantifying the impact of gravity on the preferential flow directions, then we require a flow field generated uniquely due to this effect. The problem is approached with a very similar analysis to the one presented in section \ref{Sec curvature dissipation C2 general}. An appropriate variation on the flow hypothesis and the same construction of finite volume mesh for the triangulation of a piecewise $C^{1}$ surface will be used.
%
%
%
%
\subsection{Flow Hypothesis}\label{Sec gravity flow hypothesis}
For a given triangulation $\Gamma_{\scriptscriptstyle \Hp, \,\triang}$ of the surface $\Gamma$, the idealized flow field must experience changes only due to the relative difference of heights of the flat faces of the elements of the triangulation. Hence, it must satisfy the following conditions
\begin{enumerate}[(i)]
\item
The velocity is constant in magnitude and direction within a flat face.

\item
The magnitude of the velocity within a flat element $\Kz$ is given by
\begin{subequations}\label{Def hypothesis of velocity due to gravity}
\begin{equation}\label{Def magnitude of velocity due to gravity}
\left\vert \u(K)\right\vert\defining\sqrt{2\,g\,\left(\frac{1}{2\,g}+\zmax - \pk\cdot\kversor\right)}.
\end{equation}
Where $g$ is the gravity and
\begin{equation}\label{Def maximum height}
\zmax\defining\max\left\{\pl\cdot\kversor:L\in \triang\right\} .
\end{equation}
\end{subequations}
The addition of the quantity $\frac{1}{2\,g}$ to the height of reference in the expression \eqref{Def magnitude of velocity due to gravity} above is meant to have velocities of magnitude $1$ on the volumes of control $\Lz$ of highest altitude (where the velocity has minimum magnitude).

\item
The field must meet the continuity flow condition, therefore the discharge rate $Q_{0}$ must remain constant. Thus, for all $K, L\in \triang$ it must hold
\begin{equation}\label{Def gravity continuity condition}
w(K)\left\vert \u(K)\right\vert =  w(L)\left\vert\u(L)\right\vert \equiv  Q_{0}  \defining 1 .
\end{equation}
Here $w(K)$ indicates the height of the fluid layer on the element $\Kz$. For simplicity we defined the discharge rate to be $1$. 

\end{enumerate}
\begin{remark}\label{Rem the height variable}
   Observe that in condition (iii) above we needed to introduce the fluid layer height $w(L)$. However it is the reciprocal of the velocity magnitude $\vert \u(L)\vert$ i.e. it is not and independent variable. 
\end{remark}
%
%
%
%
\subsection{The Velocity Field Due to Gravity}\label{Sec gravity field flow model}
Using the construction of the velocity field presented in section \ref{Sec velocity field construction} we have that for a given \emph{master velocity} $\ugamma = \alpha_{\,1}\,\iversor+\alpha_{\,2}\,\jversor$ and $K\in \triang$, the velocity at the volume of control $\Kz$ is given by
\begin{subequations}\label{Eq velocity due to gravity}
\begin{equation}\label{Eq local velocity due to gravity}
\u(K)
\equiv \sqrt{2\,g\left(\frac{1}{2\,g}+\zmax - \pl\cdot\kversor\right)}
\left(\begin{array}{cc}
1- \dfrac{\n_{1}^{\,2}}{1+\n_{\,3}}   & -\dfrac{\n_{1}\,\n_{\,2}}{1+\n_{\,3}}\\[8pt]
-\dfrac{\n_{1}\,\n_{\,2}}{1+\n_{\,3}} & 1- \dfrac{\n_{\,2}^{\,2}}{1+\n_{\,3}}\\[8pt]
\n_{1} & \n_{\,2}
\end{array}
\right)
\left\{\begin{array}{c}
\alpha_{\,1}\\[18pt]
\alpha_{\,2}
\end{array}
\right\} 
\defining \left(\mathcal{G}_{\Gamma_{\Hp,\,\triang}}\ugamma \right) (K) .
\end{equation}
With $\n(K) = \n_{1} \iversor +  \n_{2}\jversor + \n_{3} \kversor$ . The global velocity field $\mathcal{G}_{ \Gamma_{ \Hp,\,\triang}}\ugamma : \Hp\rightarrow \R^{\!3}$ is defined by
\begin{equation}\label{Eq global velocity due gravity 3-D}
 \mathcal{G}_{ \Gamma_{ \Hp,\,\triang}}\ugamma(s) \defining
 \sum_{K\in \triang} \u(K) \, \ind_{K}(s)\,,\quad s\in \Hp .
\end{equation}
\end{subequations}
Again, the average operator has analogous properties
\begin{lemma}\label{Th average velocity operator on gravity field}
   Let $\Gamma$ be a piecewise $C^{1}$ surface, $\Gamma_{\scriptscriptstyle \Hp, \triang}$ be a triangulation and the average velocity operator $m$ defined by \eqref{Def average velocity operator}, then 
   \begin{enumerate}[(i)]
   \item \label{Th null kernell average velocity operator on gravity field}
   $\ker(m\circ  \mathcal{G}_{\Gamma_{\Hp,\,\triang}}) = \{\zero\}$. 
   
   \item \label{Th rank average velocity operator on gravity field}
   The space $(m\circ  \mathcal{G}_{\Gamma_{\Hp,\,\triang}})(\R^{\!2})$ is two dimensional.  
   \end{enumerate}
   \begin{proof}
   Identical to the proof of lemma \ref{Th average velocity operator}.
   \end{proof}
\end{lemma}
%
%
%
\subsection{Dissipation of Mechanical Energy Due to Gravity}\label{Sec loss due to gravity velocities}
We compute the loss of mechanical energy due to gravity applying the procedure presented in section \ref{Sec dissipation due to curvature}, but using the flow field defined by the equations \eqref{Eq velocity due to gravity}. It is important to observe that the outcome is not a multiple of the previous case given by the tensor $D(\sigma)$ in \eqref{Def discreet strain rate tensor}, because the difference $\u(K) - \u(L)$ has new values of velocity magnitude, although the geometric characteristics of the triangulation are preserved.
%
%
%
\subsection{The Related Eigenvalue Problem}\label{Sec related eigenvalue problem}
The analysis made in section \ref{Sec mininimal dissipation due to curvature} is entirely applicable in the current case i.e. for a master velocity $\b\in \R^{\,2}$, the dissipation of mechanical energy due to the induced velocity field is given by
\begin{equation*}
U_{grav} (\b)= \b^{T} M_{grav}\,\b
\end{equation*}
With $M_{grav}$ symmetric, positive semi-definite matrix. Again, the spectral theorem yields the existence of an orthonormal basis of eigenvectors $\{\fversor_{\!\!\! 1}, \fversor_{\!\!2}\}$ such that
\begin{equation}\label{Eq gravity internal energy}
U_{grav} (\b)= \lambda_{1}\left(\b\cdot\fversor_{\!\!\!1}\right)^{2}+\lambda_{\,2}\left(\b\cdot\fversor_{\!\!2}\right)^{2} .
\end{equation}
For $0\leq \lambda_{1}\leq \lambda_{2}$. Define the minimizing set
\begin{equation}\label{Eq minimizing set due to gravity}
\omega_{grav}\defining\left\{\b\in S^{1}:U_{grav}(\b) = \lambda_{1}\right\} .
\end{equation}
The set of Preferential Fluid Flow Directions of the surface $\Gamma_{\Hp, \triang}$ \emph{due to gravity} has to be contained in
\begin{equation}\label{Eq potential preferential direction due to gravity}
\left\{\frac{m(\mathcal{G}_{\Gamma_{\Hp,\,\triang} }\b)}{\vert m(\mathcal{G}_{\Gamma_{\Hp,\,\triang} }\b) \vert}:\b\in \omega_{grav}\right\} .
\end{equation}
As in section \ref{Sec preferential direction due to curvature} there are two possible cases. If $\lambda_{1}<\lambda_{2}$, then $\omega_{grav}$ has only two points, namely $\fversor_{\!\!\!1}$ and $-\fversor_{\!\!\!1}$ for $\fversor_{\!\!\! 1}$ the unitary eigenvector associated to $\lambda_{1}$. If $\lambda_{1} = \lambda_{2}$, then the set $\omega_{grav}$ has infinitely many elements. However, in both cases a further analysis has to be made in order to determine the preferential flow directions.
%
%
%
%
\subsection{The External Mechanical Energy of the Fluid and the Entropy Choice Function}\label{Sec external mechanical energy and entropy}
%
%
%
So far $U_{grav}(\b)$ accounts for the total \emph{internal} energy dissipation, and due to \eqref{Eq gravity internal energy} $U_{grav}(\b) = U_{grav}(-\b)$ for all $\b\in \R^{2}$, therefore a final criterion needs to be set in order to decide which direction in $\{\b, -\b\}$ is preferential over the other, or if none of them is. First we introduce a definition
\begin{definition}\label{Def outer normal of triangles}
For any $\Kz\subset \Gamma_{\scriptscriptstyle \Hp, \,\triang}$ we denote $\{\nukl:1\leq \ell\leq 3\}\subset \R^{\,3}$ the outward normal vectors of the edges such that $\nukl\cdot\n(K) = 0$.
If $\sigma$ is an edge of the element $K$ we will use $\nusigk$ to designate the outer normal vector perpendicular to $\sigma$ and $\n(K)$.
\end{definition}
The total external energy of the free fluid is given by the algebraic sum of kinetic and potential energies, this is
\begin{multline}\label{Eq total external energy}
E_{grav} (\ugamma) 
\defining\frac{1}{2}\, \rho \sum_{K\,\in\, \triang}\vert\Kz\vert \, w(K)
\sum_{\substack{\sigma\, \in \E, \, \sigma \subseteq \partial K\\ \nusigk\cdot\,\u(K)>\,0}}
\frac{\nusigk\cdot\,\u(K) }
{\sum
\{\nutauk\cdot\,\u(K):\tau\in\partial K, \,\nutauk\cdot\,\u(K)>\,0\}}  \, \left\vert \nusigk\cdot\,\u(K) \right\vert^{2} \\
+\rho\, g\sum_{K\,\in\, \triang}\vert\Kz \vert \,  w(K)
\sum_{\substack{\sigma\, =\, K\vert L\\\nusigk\cdot\,\u(K)>\,0}}
\frac{\nusigk\cdot\,\u(K) }
{\sum
\{\nutauk\cdot\,\u(K):\tau\in\partial K, \,\nutauk\cdot\,\u(K)>\,0\}} 
\left(\pl - \pk\right)\cdot\kversor \\
+\rho\, g\sum_{K\,\in\, \triang}\vert\Kz  \vert \, w(K)
\sum_{\substack{\sigma\, \in \E - \E_{int}\\ \sigma \subseteq \partial K, \, \nusigk\cdot\,\u(K)>\,0}}
\frac{\nusigk\cdot\,\u(K) }
{\sum
\{\nutauk\cdot\,\u(K):\tau\in\partial K, \,\nutauk\cdot\,\u(K)>\,0\}} 
\left(\qs - \pk\right)\cdot\kversor .
\end{multline}
In the expression above $g$ is the gravity, $\rho$ the fluid density and $\left\{\u(K):K\in \triang\right\}$ is the flow field defined by \eqref{Eq local velocity due to gravity}. The approximations for kinetic and potential energy are given by $\frac{1}{2}\, \rho \, \vert\Kz  \vert \,  w(K)\vert \nusigk\cdot \u(K)\vert^{2}$ and $\rho \, g \,\vert\Kz \vert\,  w(K) \left(\pl - \pk\right)\cdot\kversor$ respectively. The condition $\nusigk\cdot\u(K)>0$ chooses which edges of the element $K$ are ``downstream''. The kinetic energy is quantified only in the summand of the first line of the right hand side in \eqref{Eq total external energy}. However, the potential energy needs to be quantified using sub-cases, in \eqref{Eq total external energy} the second line accounts for the interior edges while the third line for the exterior edges. The latter uses the transmission point $\qs$ defined in \eqref{Def triangulation points of transmission}. All the quantities are weighted by the factor:
\begin{equation}\label{Eq mass distribution factor}
\left\{\sum_{\substack{\tau\,\in\, \partial K\\\nutauk\cdot\,\u(K)>\,0}}\nutauk\cdot\,\u(K)\right\}^{-1}\nusigk\cdot\,\u(K).
\end{equation}
The weight above accounts for the fraction of fluid mass flowing through the edge $\sigma$.
\begin{remark}\label{Rem entropy choice function}
Observe that unlike $U_{grav}$ the function $E_{grav}(\cdot)$ is not even because of the ``downstream choice" $\nusigk\cdot\u(K)>0$. More precisely since $\#\{\sigma \subseteq K: \nusigk\cdot\u(K)>0\}\in \{1, 2\}$ we have
\begin{equation*}
\begin{split}
   \text{If}\;\;\#\{\sigma \subseteq K: \nusigk\cdot(\u(K) )>0\}  = 2 \, ,
   \quad \text {then} \;\; \#\{\sigma \subseteq K: \nusigk\cdot(-\u(K) )>0\}  = 1\, .\\
   \text{If}\;\;\#\{\sigma \subseteq K: \nusigk\cdot(\u(K) )>0\}  = 1 \, ,
   \quad \text {then} \;\; \#\{\sigma \subseteq K: \nusigk\cdot(-\u(K) )>0\}  = 2 .
   \end{split}
\end{equation*}
 Consequently $E_{grav}(- \ugamma)\notin\{ E_{grav}(\ugamma),
- E_{grav}(\ugamma) \}$, i.e. $E_{grav}(\cdot)$ is neither even, nor odd.
\end{remark}
Finally, we use $E_{grav}$ in \eqref{Eq total external energy} to decide the preferential direction
\begin{definition}\label{Def entropy choice definition}
Let $\echoice:\R^{\,2}\rightarrow \R^{\,2}$ be the \emph{Entropy Choice Function} defined as follows
\begin{equation}\label{Def entropy choice function}
\echoice(\b) =  \begin{cases}
\b & E_{grav}(\b) = \max\left\{E_{grav}(\b), E_{grav}(-\b)\right\} ', ,\\
-\b & E_{grav}(\b) < \max\left\{E_{grav}(\b), E_{grav}(-\b)\right\}\, .
\end{cases}
\end{equation}
%
\end{definition}
Given the fact that the flow fields induced by $\b$ and $-\b$ satisfy $U_{grav}(\b) = U_{grav}(-\b)$ as seen in \eqref{Eq gravity internal energy}, the \emph{Entropy Choice Function} states the flow occurs \emph{Preferentially} on one direction over the other. This way, it provides the choice when the \emph{Internal Energy States} of the \emph{flow configurations} are identical \cite{SearsSalinger}.
%
%
%
%
\subsection{Preferential Fluid Flow Directions due to Gravity and Their Probability Space}\label{Sec preferential directions due to gravity}
%
%
%
Having $\echoice$ at our disposal we define the preferential directions by 
\begin{equation}\label{Eq preferential directions due to gravity}
\varpi_{grav}\defining\left\{\frac{m( \mathcal{G}_{\Hp, \triang}(\b) )}
{\vert m( \mathcal{G}_{\Hp, \triang}(\b) )\vert}:\b\in \omega_{grav}
,\,\b = \echoice (\b)\right\} .
\end{equation}
As for the probabilistic distribution we have three possible cases.

\emph{Case 1.} $\lambda_{1} < \lambda_{\,2}$, $\fversor_{\!\!\! 1} \neq \echoice(\fversor_{\!\!\! 1})$ or $-\fversor_{\!\!\!1} \neq \echoice(-\fversor_{\!\!\! 1})$. In this case the set $\omega_{grav}$ has two points and we can decide between $\fversor_{\!\!\! 1}$ or $-\fversor_{\!\!\! 1}$, without loss of generality we assume $\fversor_{\!\!\!1} = \echoice(\fversor_{\!\!\! 1})$ then set $\varpi_{grav}$ has only one point and it has \emph{Probability} 1.
\begin{equation}\label{Def probabilistic distribution discreet gravity 1}
\prob_{grav}\left\{\frac{m( \mathcal{G}_{\Hp, \triang}(\fversor_{\!\!\! 1}) )}{\vert m( \mathcal{G}_{\Hp, \triang}(\fversor_{\!\!\! 1}) )\vert}\right\} = 1 \, .
\end{equation}

\emph{Case 2.} $\lambda_{1} < \lambda_{\,2}$,  $\fversor_{\!\!\! 1} = \echoice(\fversor_{\!\!\!1})$ and $-\fversor_{\!\!\! 1} = \echoice(-\fversor_{\!\!\! 1})$. In this case the set $\varpi_{grav}$ has exactly two points and we can not decide between them. Therefore we adopt the \emph{Uniform Probabilistic Distribution} i.e.
\begin{equation}\label{Def probabilistic distribution discreet gravity 2}
\prob_{grav}\left\{\frac{m( \mathcal{G}_{\Hp, \triang}(\fversor_{\!\!\! 1}) )}
{\vert m( \mathcal{G}_{\Hp, \triang}(\fversor_{\!\!\! 1}) )\vert}\right\} = \frac{1}{2}\,,\quad
\prob_{grav}\left\{\frac{m( \mathcal{G}_{\Hp, \triang}(-\fversor_{\!\!\! 1}) )}
{\vert m( \mathcal{G}_{\Hp, \triang}(-\fversor_{\!\!\! 1}) )\vert}\right\} = \frac{1}{2}\, .
\end{equation}

\emph{Case 3.} $\lambda_{1} = \lambda_{\,2}$. In this case every $\eversor\in S^{\,1}$ minimizes the dissipation of internal mechanical energy of the fluid due to gravity $U_{grav}$, however $\echoice(\cdot)$ remains as a criterion to discriminate. In order to define a probability on $\varpi_{grav}$, a couple of previous results are needed
\begin{lemma}\label{Th continuity of Energy Gravitatory Function}
The function $E_{grav}: S^{1}\rightarrow \R$ is continuous.
\newline
\newline
\begin{proof}
Choose $\Kz\in \Gamma_{\Hp, \triang}$, since the map $\u\mapsto \u(K) $ is continuous the following maps are also continuous  
\begin{align*}
& \u\mapsto \u(K) \cdot \nu_{\sigma}^{K} \,,\\
 & \u  \mapsto  \nu_{\sigma}^{K} \cdot \u \;   \ind_{(0, \infty)}  (\nusigk\cdot\,\u(K) ) \, , \\[5pt]
& \u  \mapsto \nu_{\sigma}^{K} \cdot \u \;  \vert  \nu_{\sigma}^{K} \cdot \u\vert^{2} \; \ind_{(0, \infty)} (\nusigk\cdot\,\u(K) ) \, , \\[5pt]
& \u\mapsto \sum_{\substack{\sigma \in \E, \sigma \subseteq \partial K\\\nu_{\sigma}^{K} \cdot \u >0} } \nu_{\sigma}^{K} \cdot \u 
= \sum_{\sigma \in \E, \sigma \subseteq \partial K} \nu_{\sigma}^{K} \cdot \u \; \ind_{(0, \infty)} (\nusigk\cdot\,\u(K) ) .
\end{align*}
Moreover the last map above is bounded away from zero, therefore the quotient of the second over the fourth, as well as the quotient of the third over the fourth are continuous. Finally, since $E_{grav}$ is a linear combination of these type of quotients it is also continuous.
\end{proof}
\end{lemma}
\begin{theorem}\label{Th measure of eligible directions}
Consider the set 
\begin{equation}\label{Def eligible master velocities gravity}
R_{grav}\defining\left\{\eversor\in \omega_{grav}: \eversor = \echoice(\eversor) \right\} .
\end{equation}
\begin{enumerate}[(i)]

\item
The set $R_{grav}$ is measurable.

\item
The set $R_{grav}$ has positive arc-length measure in $S^{\,1}$.
\end{enumerate}
\begin{proof}
\begin{enumerate}[(i)]
\item
First observe that the function $\phi (\b) \defining \max\{E_{grav}(\b), E_{grav} (\b)\} $ is continuous and the following equality holds
\begin{equation*}
\echoice(\b) =  \b \, \left(\ind_{\{0\}} \circ( \phi - E_{grav}) \right) (\b) - \b \, \left(\ind_{(0, \infty)} \circ( \phi - E_{grav}) \right) (\b)
\end{equation*}
Then $\echoice(\cdot)$ is a measurable function since it is the composition of measurable functions. On the other hand, denoting $id$ the identity function on $S^{1}$, the set $R_{grav} $ can be written as
\begin{equation*}
R_{grav} = (id - \echoice)^{-1} \{0\}
\end{equation*}
Since $R_{grav}$ is the inverse image of a Borel set through a measurable function, it is measurable.

\item
If $\left\{\eversor\in S^{\,1}: \eversor = \echoice(\eversor) \right\} = S^{\,1}$ there is nothing to prove. If not, there exists an element in $S^{\,1}$, namely $-\eversor_{0}$ such that $-\eversor_{0}\neq \echoice(-\eversor_{0})$ i.e. $E_{grav}(-\eversor_{0})<E_{grav}(\eversor_{0})$. Notice that $F(\b)\defining E_{grav}(\b) - E_{grav}(-\b)$ is a continuous function, hence, defining $\epsilon\defining F(\eversor_{0})>0$ there must exist $\delta>0$ such that $\left\vert\eversor-\eversor_{0}\right\vert<\delta$ implies $\left\vert  F(\eversor)- F(\eversor_{0})\right\vert <\frac{\epsilon}{2}$. This implies $\frac{\epsilon}{2}< F(\eversor)$ which gives the inclusion
\begin{equation*}
\{\eversor\in S^{\,1}: \left\vert\eversor -\eversor_{0}\right\vert<\delta\}\subseteq R_{grav} .
\end{equation*}
Since an open non-empty neighborhood of $S^{\,1}$ has positive arc-length measure the result follows.
\end{enumerate}
\end{proof}
\end{theorem}
Now we endow the set $\varpi_{grav}$ with a natural $\sigma$-algebra as in definition \ref{Def sigma algebra curvature}.
\begin{definition}\label{Def sigma algebra gravity}
Let $\Gamma$ be a piecewise $C^{1}$ surface, $\Gamma_{\scriptscriptstyle \Hp, \triang}$ be a triangulation and $R_{grav}$ defined in \eqref{Def eligible master velocities gravity}. Consider the surjective function
\begin{equation}\label{Def probabilistic distribution function gravity}
\varphi:R_{grav}\rightarrow \varpi_{grav}\,,\quad
\varphi(\eversor)\defining \frac{m( \mathcal{G}_{\Hp, \triang}(\eversor) )}{\vert m( \mathcal{G}_{\Hp, \triang}(\eversor) )\vert} \, .
\end{equation}
Let $\beta$ be the family of all Borel sets of $S^{\,1}$ intersected with the set $R_{grav}$, we define the following $\sigma$-algebra
\begin{equation}\label{Eq sigma algebra gravity}
\tau_{grav}\defining\left\{A\in \wp(\varpi_{grav}): \varphi^{-1}(A)\in \beta\right\}  ,
\end{equation}
where $\wp(\varpi_{grav})$ is the power set of $\varpi_{grav}$.
\end{definition}
Finally, the \emph{Probabilistic Distribution} for the non-discrete case follows in a natural way
\begin{definition}\label{Def probabilistic distribution continuous gravity}
Consider the measure space $\left(\varpi_{grav}, \tau_{grav}\right)$, for any $A\in \tau_{grav}$ define
\begin{equation}\label{Def probabilistic measure continuous gravity}
\prob_{grav}\left[A\right]\defining \frac{1}{\mu (R_{grav})}\,\mu\left[\varphi^{-1}(A)\right].
\end{equation}
With $\mu$ the arc-length measure in $S^{\,1}$.
\end{definition}
%
%
%
%
%
%
%
%
%
%
%
%
\section{Preferential Fluid Flow due to Friction}\label{Sec preferential direction due to friction}
%
%
%
In the present section we address the question of finding a fluid flow configuration that minimizes the mechanical energy dissipation due to friction. To that end we use the same setting: $\Gamma$ a piecewise $C^{1}$ surface and  a triangulation $\Gamma_{\Hp,\triang}$. We also adopt the same flow hypothesis assumed in section \ref{Sec flow hypothesis} leading to the velocity field constructed in section \ref{Sec velocity field construction}, equation \eqref{Eq velocity due curvature 3-D}. Hence, the velocity $\u(K)$ within a volume of control $\Kz$ is constant and depends linearly with respect to a master velocity $\ugamma\in \R^{2}$.
%
%
%
\subsection{Mechanical Energy Dissipation due to Friction}\label{Sec dissipation due to friction}
The dissipation of mechanical energy due to friction is proportional to the magnitude of the velocity and to the length of the path that a fluid particle must travel i.e. it is proportional to the \emph{Length} of the \emph{Stream Lines}. 
%
%
%
\begin{figure}[!]
\caption[17]{Stream Lines, Volume of Control $\Kz$, Horizontal View}\label{Fig Stream Lines}
\includegraphics{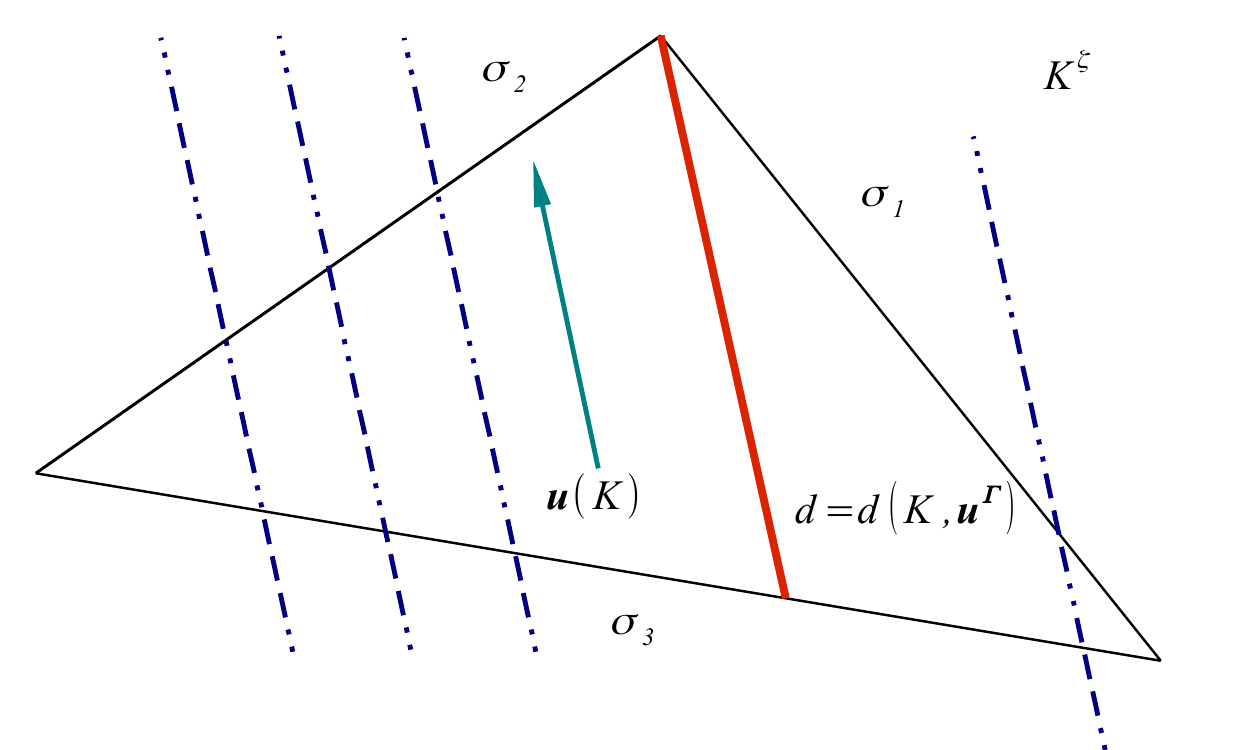}
\end{figure}
%
%
%
Figure \ref{Fig Stream Lines} depicts a fixed volume of control $\Kz$, the constant velocity $\u(K)$ within and in segmented blue lines the \emph{Stream Lines}. A simple calculation shows that the \emph{Average Length} of the \emph{Stream Lines} is given by $\frac{1}{2}\,d$ where $d = d(K, \ugamma)$ is the maximum length of the segments parallel to $\u(K)$ contained inside $\Kz$. Therefore the energy dissipation due to friction is given by
\begin{equation}\label{Eq friction dissipation functional}
\F (\ugamma)\defining \gamma \; \sum_{K\,\in\,\triang} \frac{1}{2}\,\vert \Kz \vert \,  d\left(K, \ugamma\right) \vert\u(K) \vert =
\frac{1}{2}\, \gamma \; \vert\ugamma \vert\sum_{K\,\in\,\triang}  d\left(K, \ugamma\right) \vert \Kz \vert.
\end{equation}
In the expression above $\gamma$ is the friction coefficient depending on the material of the walls of the fissure, $\vert \Kz\vert$ is the area of the volume of control $\Kz$ and $ \vert\ugamma \vert$ is the magnitude of the master velocity $\ugamma$. The second equality in \eqref{Eq friction dissipation functional} uses the fact that, by construction $\vert \u(K) \vert = \vert \ugamma\vert$ for all $K\in \triang$. The expression \eqref{Eq friction dissipation functional} implies that we need to address the question
\begin{equation}\label{Eq minimization problem friction}
\vgamma\in S^{\,1}:
\F(\vgamma) = \min\left\{\F(\b): \b \in S^{\,1}\right\} .
\end{equation}
However, the functional $\F\left(\cdot\right)$ does not have the same quadratic structure as $U_{curv}$ or $U_{grav}$ due to the presence of the map $\ugamma\mapsto d(K, \ugamma)$. Although proving the existence of minima is easy to show, characterizing them is extremely difficult. Nevertheless it will be shown that the set of preferential directions can be turned into a suitable probability space. To that end, the function $\ugamma\mapsto d(K, \ugamma)$ needs to be further studied.
%
%
%
\subsection{The Minimizing Set of $\F$}\label{Sec Convexity of Inner Segment Length Function}
\begin{figure}[!]
\caption[18]{Triangle $\&$ Inner Segments}\label{Fig Triangle and Inner Segments}
\includegraphics{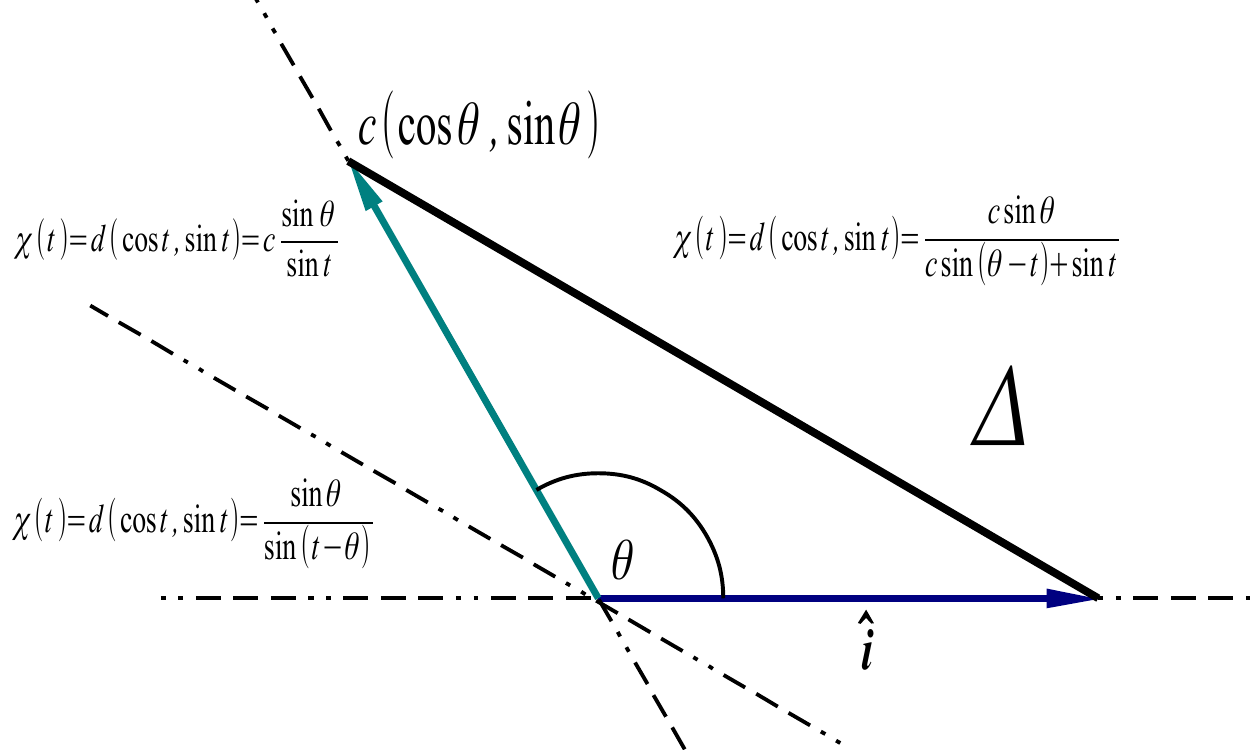}
\end{figure}

We start analyzing a problem closely related to the map $\ugamma\mapsto d(K, \ugamma)$, introducing an auxiliary function 
\begin{definition}\label{Def parametrized inner segment length function}
   The parametrized inner-segment length function $\chi: \R\rightarrow [0, \infty)$ is defined by  $$\chi(t) \defining d(\cos t, \sin t).$$
\end{definition}
It is immediate to see that $\chi$ is $\pi$-periodic. Now consider the following family of particular triangles in $\R^{2}$, see figure \ref{Fig Triangle and Inner Segments}.
\begin{multline}\label{Def flat quasi general triangles in 2-D}
\mathcal{C} \defining \{\Delta\subseteq \R^{\!2} : \Delta \; \, \text{is a triangle with one of its sides defined by the unitary vector} \; \, \iversor\\
\text{ and a second side defined by a vector} \; c\, (\cos \theta, \sin \theta) \text{ for } 0< \theta < \pi , \, c > 0\} .
\end{multline}
With the definition above we have the following result
\begin{proposition}\label{Th  inner segment length in work triangles}
Let $\Delta$ be an element of $\mathcal{C}$ 
\begin{enumerate}[(i))] 
   \item 
The parametrized inner-segment length function $\chi: [0, \pi]\rightarrow \R$ is described by the following expression
\begin{equation}\label{Eq inner segment length in work triangles}
\chi (t) = \frac{c \, \sin \theta }{c \, \sin (\theta - t)  + \sin t}\, \ind_{[0, \theta] } 
+ c\, \frac{ \sin \theta}{\sin t }\, \ind_{[\theta, \, t_{0} ] } 
+ \frac{\sin \theta}{\sin (t - \theta)} \ind_{[t_{0}, \, \pi ] } .
\end{equation}
Where $t_{0} \in [0, \pi]$ is such that $(\cos t_{0}, \sin t_{0})$ is parallel to the third side $  (c\, \cos \theta - 1, c \, \sin \theta )$.

\item
The function $\chi$ is continuous.
\end{enumerate}
\begin{proof}
(i) By direct calculation. (ii) Follows directly from the expression \eqref{Eq inner segment length in work triangles} in the previous part.
\end{proof}
\end{proposition}
\begin{remark}\label {Rem inner segment length in work triangles}
   Observe that any triangle in $\Delta \subseteq \R^{2}$ can be reduced to an element of $\mathcal {C}$ by dilation and rotation maps. The first map only amplifies the values of $\chi$, while the second map only shifts the sub-intervals of definition, which will become at most four. In all the cases the function $\chi$ is described by rational trigonometric functions bounded away from singularities. 
\end{remark}
In order to extend the analysis to the elements of the triangulation we define a new family
\begin{multline}\label{Def flat quasi general triangles in 3-D}
\mathcal{D} \defining \{\Delta\subseteq \R^{\!3} : \Delta \; \, \text{is a triangle whose orthogonal projection onto the horizontal plane}\,\; \langle \kversor \rangle^{\perp} , \\
%
\text{is an element of} \;\, \mathcal {A}\} .
\end{multline}
Now we describe the function $\chi$ on the family $\mathcal{D}$. 
\begin{proposition}\label{Th  inner segment length in work triangles 3 - D}
Let $\Delta$ be an element of $\mathcal{D}$, $\n $ be the normal vector orthogonal to $\Delta$ and let $\Delta_{0}$ be the orthogonal projection of $\Delta$ onto the horizontal plane $\langle \kversor \rangle ^{\perp}$. 
\begin{enumerate}[(i)]
\item
Let $\chi, \chi_{0}: [0, \pi]\rightarrow \R$ be the respective parametrized inner-segment length functions then, the following relation holds
\begin{equation}\label{Eq inner segment length in work triangles 3 - D}
  \chi(t)  = \frac{1}{\sqrt{1 - (\n \cdot (\cos t, \sin t, 0) )^{2}  } }  \,\chi_{0} (t) .
\end{equation}

\item
The function $\chi:[0, \pi]\rightarrow [0, \infty)$ is continuous.
\end{enumerate}
\begin{proof}
\begin{enumerate}[(i)]
\item 
It suffices to observe that the factor affecting $\chi_{0} (t)$ in \eqref{Eq inner segment length in work triangles 3 - D} is the amplification factor on the segment length (or vector norm) due to the ``lifting'' of the segment $\chi_{0}(t)(\cos t, \sin t)$ contained in $\Delta_{0}$ onto the triangle $\Delta$.

\item
Follows immediately from the continuity of $\chi_{0}$ shown in part (ii) of proposition \ref{Th  inner segment length in work triangles} and the relation \eqref{Eq inner segment length in work triangles 3 - D} given by the previous part.

\end{enumerate}
\end{proof}
\end{proposition}
\begin{remark}\label {Rem inner segment length in work triangles 3-D}
   Observe that any triangle in $\Delta \subseteq \R^{3}$ can be reduced to an element of $\mathcal {B}$. Denote $\Delta_{0}$ the orthogonal projection of $\Delta$ on the horizontal plane $\langle \kversor \rangle^{\perp}$. Then, the dilation and rotation maps that turn $\Delta_{0}$ into an element of $\mathcal{C}$ also turn $\Delta$ into an element of $\mathcal{D}$. 
\end{remark}
%
%
%
%
%
%
\begin{theorem}\label{Th existence of the minimum frictional functional}
   Let $\Gamma$ be a piecewise $C^{1}$ surface and $\Gamma_{\scriptscriptstyle \Hp, \triang}$ be a triangulation. Then, the functional $\F$ attains its minimum on $S^{1}$.
   \newline
   \newline
   \begin{proof}
      Using the definition of $\F$, consider the function $X: [0, \pi]\rightarrow [0, \infty)$
\begin{equation*}
X(t) \defining \F (\cos t, \sin t)  =
\frac{1}{2}\, \gamma \sum_{K\,\in\,\triang}  d\left(K, (\cos t, \sin t)\right) \vert \Kz \vert.
\end{equation*}
Since $X$ is the sum of continuous functions as seen in proposition \ref{Th  inner segment length in work triangles 3 - D} part (ii), the function $X$ is continuous. Since it is defined on the compact set $[0, \pi]$ it must attain its minimum at a point, namely $t_{0}$. Finally, since
\begin{equation*}
X(t_{0}) = \min \{X(t): t\in [0, \pi]\} \equiv \min  \{\F (\ugamma)  : \ugamma\in S^{1}\}
= \F(\cos t_{0}, \sin t_{0}) ,
\end{equation*}
the result follows.
   \end{proof}
\end{theorem}
\begin{lemma}\label{Th analytic extenion existence}
Let $\Delta\subseteq \R^{\!3}$ be a triangle, $\chi:[0, \,\pi]\rightarrow [0, \infty)$ its parametrized inner-segment length function and $\{I_{n}: 1\leq n\leq N \}$ the disjoint open sub-intervals of piecewise definition, with $N\in \{3, 4\}$. Then for each $I_{n}$ the function $\chi\vert_{ I_{n}}$ has an analytic extension to an open connected set $U_{n}$ of the complex plane, containing the closure of $I_{n}$.
\newline
\newline
\begin{proof}
   Due to propositions \ref{Th  inner segment length in work triangles} and \ref{Th  inner segment length in work triangles 3 - D} as well as remarks \ref{Rem inner segment length in work triangles}, \ref{Rem inner segment length in work triangles 3-D} we observe that the function $\chi\vert_{ I_{n}}$ is a rational trigonometric function bounded away from singularities. Therefore, the map $z\mapsto \chi(z)$ is well-defined and analytic in a neighborhood of the closure of $I_{n}$ for each $1\leq n \leq N$.
\end{proof}
\end{lemma}
Finally we present the result that will allow us to define the uniform probability on the set of preferential directions. 
\begin{theorem}\label{Th finite or non trivial measurable directions for friction}
   Let $\Gamma$ be a piecewise $C^{1}$ surface and $\Gamma_{\scriptscriptstyle \Hp, \triang}$ be a triangulation. Then, the minimization set 
\begin{equation}\label{Def frictional minimization set}
\omega_{friction} \defining \left\{\ugamma\in S^{1}:\F(\ugamma) = \min\left\{\F(\vgamma):\vgamma\in S^{1}\right\}\right\} ,
\end{equation}
has either a finite number or elements or, it has positive Lebesgue measure. 
\newline
\newline
\begin{proof}
Let $K$ be an element of $\Gamma_{\Hp, \triang}$ and $\chi^{K}:[0, 2\,\pi]\rightarrow [0, \infty)$ be its parametrized inner-segment length function, extended by periodicity from $[0, \pi]$ to $[0, 2\pi]$. Let $\{t_{j}^{K}: 1\leq j\leq J(K)\}$ be the points in $ (0, 2\,\pi)$ connecting the disjoint sub-intervals $\{I_{n}^{K}: 1\leq n\leq N(K)\}$ of definition for $\chi_{K}$; with $J(K) = N(K) -1$ and $N(K) \in \{5, 6\}$. Define the set
   \begin{equation*}
   D \defining \{0, \, 1\}\bigcup_{K\in \triang}\{t_{j}^{K}: 1\leq j\leq J(K)\} .   
   \end{equation*}
   This set is finite and therefore it can be ordered monotonically $0 = t_{0} < t_{1} <\ldots < t_{J} = 2\pi$. Denote $I_{n} \defining (t_{j-1}, t_{j})$ for all $1\leq j\leq J$. Let $\{U^{K}_{n}: 1\leq n \leq N(K)\}$ be the open sets of analytic extension for $\left.\chi^{K}\right\vert_{I_{n}}$ such that $cl (I^{K}_{n})\subseteq U^{K}_{n}$. Define the following sets
\begin{equation}\label{Def sets of local analytic extension for F}
U_{j}\defining \bigcap \{U^{K}_{n}: I_{j} \subseteq I^{K}_{n}, \, K\in \triang\}\, , \quad 1\leq j\leq J.
\end{equation}
 The set $U_{j}$ is open connected since it is the finite intersection of open connected sets and since $cl (I_{j})\subseteq cl (I^{K}_{n})\subseteq U^{K}_{n}$ for all the chosen open sets in \eqref{Def sets of local analytic extension for F}, this implies that $cl (I_{j})\subseteq U_{j}$. Also observe that $\{U_{j}: 1\leq j\leq J\}$ is an open covering of the interval $[0, 2\,\pi]$. Consider the function $X_{j}: U_{j}\rightarrow \C$
\begin{equation*}
X_{j}(z) \defining 
\frac{1}{2}\, \gamma \sum_{K\,\in\,\triang}  \left. \chi^{K}\right\vert_{U_{j}}(z) \; \vert \Kz \vert
- \min \{\F(\ugamma): \ugamma\in S^{1} \}\, , \quad 1 \leq j \leq J .
\end{equation*}
This function is also analytic because the map $\left. \chi^{K}\right\vert_{U_{j}}$ is analytic for all $1\leq j\leq J$. If the minimization set $\omega_{friction}$ defined in \eqref{Def frictional minimization set} contains a finite number of points there is nothing to prove. If it contains infinitely many points, then it has a convergent subsequence to a point $\xi \in [0, 2\,\pi]$. Let $U_{j_{0}}$ be one of the sets in the covering $\{U_{j}: 1\leq j\leq J\}$ such that $\xi \in U_{j_{0}}$ then, the zeros of the function $X_{j_{0}}$ have a limit point, therefore it must be constant and $\omega_{friction} \supseteq I_{j_{0}}$. Recall that $ \omega_{friction} = \F^{-1} (\{\min_{S^{1}} \F\})$ is the inverse image of a point through a continuous function and therefore it is a measurable set. Seeing that the interval $I_{j_{0}}$ is non-degenerate, we conclude that $\omega_{friction}$ has positive Lebesgue measure. 
\end{proof}
\end{theorem}
%
%
%
%
%
%
\subsection{Probability Space of Preferential Fluid Flow Directions due to Friction}\label{Sec preferential directions and distribution}
We follow the same reasoning presented in section \ref{Sec preferential direction due to curvature}. 
The \emph{Preferential Fluid Flow Directions} of the surface $\Gamma_{\Hp,\,\triang}$ due to \emph{Friction} are given by
\begin{equation}\label{Def friction preferential direction}
\varpi_{friction} \defining \left\{\frac{m(V_{\Gamma_{\Hp,\,\triang} }\b)}{\vert m(V_{\Gamma_{\Hp,\,\triang} }\b)  \vert}: \b \in \omega_{friction}\right\} .
\end{equation}
Now we define the natural $\sigma$-algebra on $\varpi_{friction}$.
\begin{definition}\label{Def sigma algebra friction}
Let $\Gamma$ be a piecewise $C^{1}$ surface, $\Gamma_{\scriptscriptstyle \Hp, \triang}$ be a triangulation and the minimizing set $\omega_{friction}$ defined in \eqref{Def frictional minimization set}. Consider the surjective function 
\begin{equation}\label{Def probabilistic distribution function friction}
\varphi:\omega_{friction} \rightarrow \varpi_{curv}\,,\quad
\varphi(\eversor)\defining \frac{m(V_{\Gamma_{\Hp,\,\triang} }\eversor)}{\vert m(V_{\Gamma_{\Hp,\,\triang} }\eversor) \vert} .
\end{equation}
Let $\beta$ be the Borel tribe in $S^{1}$ intersected with $\omega_{friction}$, define the following $\sigma$-algebra
\begin{equation}\label{Eq sigma algebra friction}
\tau_{friction}\defining
\left\{A \in \wp(\varpi_{friction}): \varphi^{-1}(A) \in \beta\right\}.
\end{equation}
Where $\wp(\varpi_{friction})$ is the power set of $\varpi_{friction}$.
\end{definition}
Finally, we endow the preferential flow space with the uniform probability distribution 
\begin{definition}\label{Def distribuition due to friction}
Let $\Gamma$ be a piecewise $C^{1}$ surface, $\Gamma_{\scriptscriptstyle \Hp, \triang}$ be a triangulation and $\varpi_{friction}$ be the associated preferential fluid flow directions defined in \eqref{Def friction preferential direction} then
\begin{enumerate}[(i)]
\item
If $\varpi_{friction}$ is finite, define
\begin{equation}\label{Def dicrete probabilistic distribution friction}
\prob_{friction}\left\{\eversor\right\} \defining \frac{1}{\# \, \varpi_{friction}}\,,\quad
\forall\; \eversor \in \varpi_{friction}.
\end{equation}

\item
If $\varpi_{friction}$ is infinite, define
\begin{equation}\label{Def continuous probabilistic distribution friction}
   \prob_{friction}[A] \defining \frac{1}{\mu (\varphi^{-1} (\varpi_{friction} ) ) }
   \, \mu (\varphi^{-1} ( A ) ) = \frac{1}{\mu (\omega_{friction}  ) }
   \, \mu (\varphi^{-1} ( A ) )\, , \quad \forall \, A \in \tau_{friction}.
\end{equation}
Here $\mu$ indicates the arc-length measure in $S^{1}$.
\end{enumerate}
\end{definition}
\begin{remark}\label{Rem well defined space of friction directions}
Due to the theorem \ref{Th finite or non trivial measurable directions for friction} the probability given in the definition above is well-defined. Hence $(\varpi_{friction}, \tau_{friction}, \prob_{friction})$ is a probability space.
\end{remark}
%
%
%
%
%
%
%
%
%
%
%
\section{The Global Space of Preferential Fluid Flow Directions and Closing Remarks}
%
%
%
In this section we use the \emph{Superposition Principle} to assemble the effect of the factors analyzed in the previous sections. 
%
%
%
%
\subsection{Global Space of Preferential Fluid Flow Directions}\label{Sec Global Space of Preferential Fluid Flow Directions}
In the following assume that $\tau_{curv}$, $\tau_{grav}$ and $\tau_{friction}$ are the $\sigma$-algebras corresponding to each preferential set $\varpi_{curv}$, $\varpi_{grav}$ and $\varpi_{friction}$. If one of the preferential sets is discrete it is assumed that the associated $\sigma$-algebra is its power set. Now we introduce some definitions.
\begin{definition}\label{Def space of global preferential directions}
Let $\Gamma$ be a piecewise $C^{1}$ surface and $\Gamma_{\scriptscriptstyle \Hp, \triang}$ be a triangulation. The minimizing space associated to the triangulation $\Gamma_{\Hp, \triang}$ is defined by
\begin{subequations}\label{Def minimizing set space}
\begin{equation}\label{Def minimizing set triangulation}
\omega\left(\Gamma_{\Hp, \triang}\right)\defining\varpi_{curv}\times\
\varpi_{grav}\times\varpi_{friction} \, ,
\end{equation}
endowed with the product $\sigma$-algebra 
\begin{equation}\label{Def minimizing sigma algebra}
\tau_{\omega} \defining \tau_{curv}\otimes \tau_{grav} \otimes \tau_{friction} \, ,
\end{equation}
and the product probability measure 
\begin{equation}\label{Def minimizing probability}
\prob_{\!\!\! \omega}\defining \prob_{\!\!\! curv} \otimes\prob_{\!\!\! grav} \otimes \prob_{\!\!\! friction} \, .
\end{equation}
\end{subequations}
\end{definition}
\begin{definition}\label{Def weights of preferential directions}
We define the weights of the preferential directions coming from each effect in the following way
\begin{subequations}\label{Eq weights of preferential directions}
\begin{equation}\label{Eq weight curvature}
p_{\,1}\defining \frac{\min U_{curv}}{\min U_{curv} + \min U_{grav} + \min \,\F} \, ,
\end{equation}
\begin{equation}\label{Eq weight gravity}
p_{\,2}\defining \frac{ U_{width}+\min U_{grav}}{ U_{curv}+ \min U_{grav} + \min \,\F} \, ,
\end{equation}
\begin{equation}\label{Eq weight friction}
p_{\,3}\defining \frac{\min \,\F}{\min U_{curv}+ \min U_{grav} + \min \,\F} \, .
\end{equation}
\end{subequations}
\end{definition}
Next we define the \emph{Space of Preferential Directions}
\begin{definition}\label{Def superposition function and preferential set}
Let $\Gamma$ be a piecewise $C^{1}$ surface, $\Gamma_{\scriptscriptstyle \Hp, \triang}$ be a triangulation and consider the function $\Phi: \varpi_{curv}\times\varpi_{grav}\times\varpi_{friction}\rightarrow S^{\,2}\cup\{\zero\}$ defined by
\begin{equation}\label{Eq superposition function}
\Phi(\a_{\,1}, \a_{\,2}, \a_{\,3})\defining 
\begin{cases}\dfrac{p_{\,1}\,\a_{\,1}+p_{\,2}\,\a_{\,2}+p_{\,3}\,\a_{\,3}}
{\left\vert p_{\,1}\,\a_{\,1}+p_{\,2}\,\a_{\,2}+p_{\,3}\,\a_{\,3}\right\vert}  & p_{\,1}\,\a_{\,1}+p_{\,2}\,\a_{\,2}+p_{\,3}\,\a_{\,3} \neq \zero ,\\[8pt]
\zero & p_{\,1}\,\a_{\,1}+p_{\,2}\,\a_{\,2}+p_{\,3}\,\a_{\,3} = \zero . %
\end{cases}
\end{equation}
The \emph{Space of Preferential Directions} is given by 
\begin{subequations}\label{Def preferential directions space}
\begin{equation}\label{Def preferential set triangulation}
\varpi\left(\Gamma_{\Hp, \triang}\right) \defining \Phi\left(
\omega\left(\Gamma_{\Hp, \triang}\right)\right) = \Phi \left(\varpi_{curv}\times\
\varpi_{grav}\times\varpi_{friction} \right) ,
\end{equation}
endowed with the $\sigma$-algebra 
\begin{equation}\label{Def preferential sigma algebra}
\tau  \left(\Gamma_{\Hp, \triang}\right)\defining 
\{ A\subseteq \varpi\left(\Gamma_{\Hp, \triang}\right) : \Phi^{-1}(A) \in \tau_{curv}\otimes \tau_{grav} \otimes \tau_{friction} \,  ,
\end{equation}
and the probability measure 
\begin{equation}\label{Def preferential probability}
\prob_{\scriptscriptstyle \Gamma_{\Hp, \triang})} [A]\defining 
\prob_{\!\!\! \omega} [\Phi^{-1} (A) ] \, .
\end{equation}
\end{subequations}
\end{definition}
\begin{remark}\label{Rem the zero in the final space}
Observe that unlike the previous cases of analysis the direction $\zero$ is actually possible. This would mean that the medium is isotropic: it is the unlikely event that the analyzed effects cancel each other, i.e. they ``average out". 
\end{remark}
%
%
%
%
\subsection{Entropy of the Preferential Fluid Flow Information}\label{Sec Entropy of the Preferential Fluid Flow Information}
Finally, appealing to Shannon's information theory \cite{Shannon} we present a definition of the \emph{Geometric Entropy} of the triangulation $\Gamma_{\Hp,\triang}$, as a measure of the amount of uncertainty in the phenomenon of preferential fluid flow. 
\begin{definition}\label{Def finite partitions}
Denote $\Sigma\left(\Gamma_{\Hp, \triang}\right)$ the family of all countable partitions of $\varpi\left(\Gamma_{\Hp, \triang}\right)$ such that each set is $\tau\left(\Gamma_{\Hp, \triang}\right)$-measurable.
\end{definition}
The concept of entropy on the elements of $\Sigma\left(\Gamma_{\Hp,\triang}\right)$ follows naturally \cite{Shannon, Khinchin}.
\begin{definition}\label{Def partition entropy}
Let $\Gamma$ be a piecewise $C^{1}$ surface, $\Gamma_{\scriptscriptstyle \Hp, \triang}$ be a triangulation, $\Sigma\left(\Gamma_{\Hp,\triang}\right)$ the family of countable measurable partitions of $\varpi\left(\Gamma_{\Hp, \triang}\right)$ and $\A\in \Sigma\left(\Gamma_{\Hp,\triang}\right)$. 
\begin{enumerate}[(i)]
\item 
The information function associated to $\A$ is given by 
\begin{equation}\label{Eq partition information function}
I_{\scriptscriptstyle \A}(x) \defining - \sum_{A\,\in\,\A}\log \prob_{\scriptscriptstyle \Gamma_{\Hp,\triang}}\left[A\,\right] \; \ind_{A}\, (x) .
\end{equation}
\item
The \emph{Geometric Entropy} associated to $\A$ is the expectation of $I_{\A}$, i.e.
\begin{equation}\label{Eq partition entropy}
\entro\left(\A\right) \defining \int_{\varpi\left(\Gamma_{\Hp, \triang}\right)} I_{\scriptscriptstyle \A} (x) d\prob_{\scriptscriptstyle \Gamma_{\Hp,\triang}} (x) = 
- \sum_{A\,\in\,\A}\,\prob_{\scriptscriptstyle \Gamma_{\Hp,\triang}}\left[A\,\right]\log \prob_{\scriptscriptstyle \Gamma_{\Hp,\triang}}\left[A\,\right] .
\end{equation}
In the last expression it is understood that $0\cdot \log 0 = 0$ in consistency with measure theory.  

\item
The \emph{Geometric Entropy} of the triangulation $\Gamma_{\Hp,\triang}$ is given by
\begin{equation}\label{Eq triangulation entropy}
\entro\left(\Gamma_{\Hp,\triang}\right) \defining \sup \{\entro(\A): \A \in \Sigma\left(\Gamma_{\Hp,\triang}\right)\}\\
= - \int_{\varpi\left(\Gamma_{\Hp, \triang}\right)}  p(x)\, \log p(x)\, dx
\end{equation}
Where the last equality holds whenever $d\prob_{\scriptscriptstyle \Gamma_{\Hp,\triang}}$ has density $p$ with respect to the Lebesgue measure $dx$ in $S^{2}\cup \{\zero\}$. In particular for such density to exist it must hold that $\prob_{\scriptscriptstyle \Gamma_{\Hp,\triang}}(\{\zero\}) = 0$.
\end{enumerate}
\end{definition}
%
%
%
%
%
\subsection{Closing Remarks and Future Work}\label{Closing Remarks}
\begin{figure}
        \centering
        \begin{subfigure}[Meshed Fissured Medium]
                {\resizebox{6.5cm}{6.5cm}
{\includegraphics{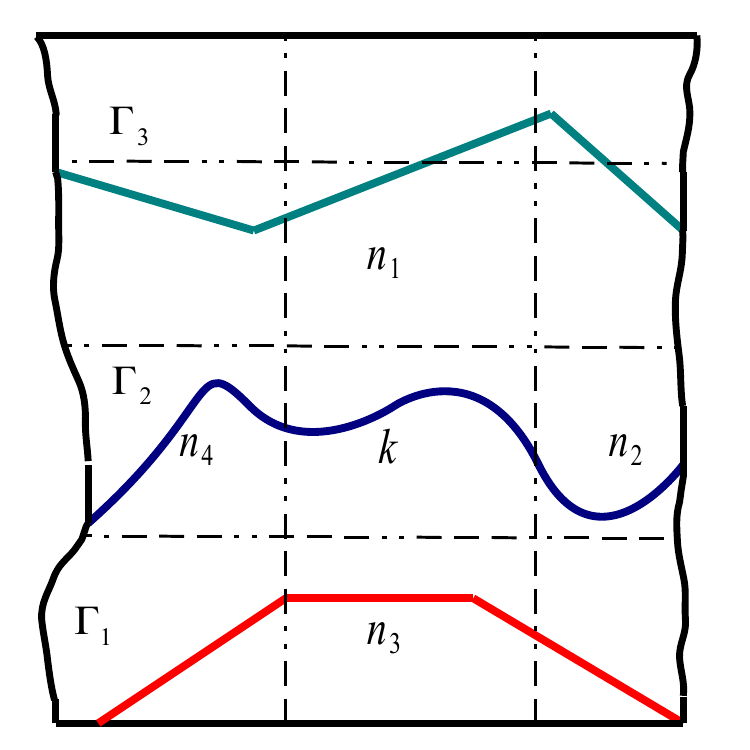} } }
                \label{Fig Meshed Geological Fissured Medium}
        \end{subfigure}
        \qquad
        ~ 
          \begin{subfigure}[Transmission Probabilities]
                {\resizebox{6.5cm}{6.5cm}
{\includegraphics{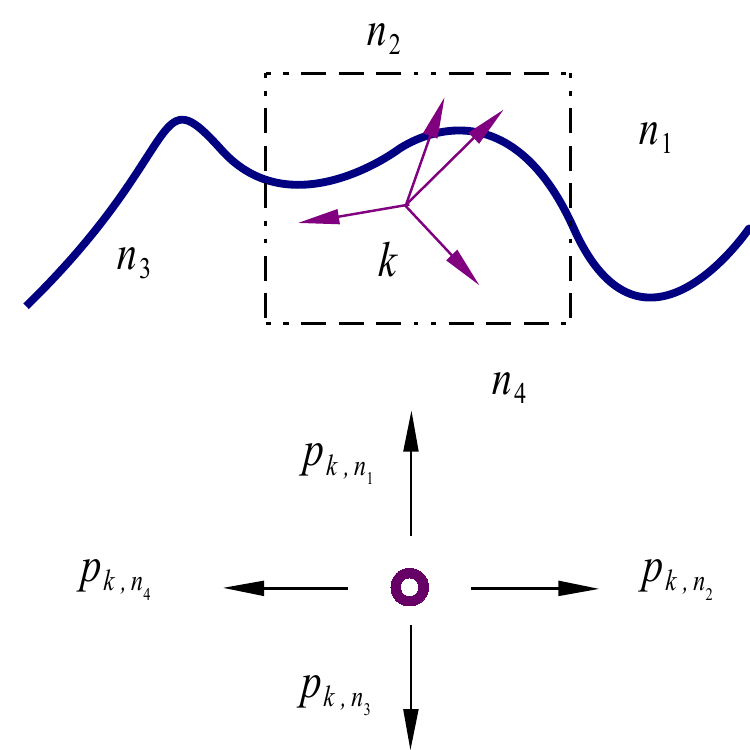} } }
\label{Fig Markov Vertex}
        \end{subfigure}%
        ~ 
\caption{Meshed Fissured Medium and Transmission Probabilities}
\end{figure}\label{Fig Transmission Probabilities}
Due to analysis exposed, several observations are in order:
\begin{enumerate}[(i)]
\item
The method we have presented to determine the probability space of preferential flow directions can be extended to other factors of interest, depending on the framework and available information.   

\item
The analyzed factors were presented in order of difficulty, not from the mathematical point of view, but from the applicability of its conclusions. Since the curvature effect is reduced to an eigenvalue-eigenvector problem this is the easiest of all. Gravity demands the inclusion of an extra criterion expensive to implement numerically, however, the minimizing set $\omega_{grav}$ is still easy to characterize by an eigenvalue-eigenvector problem. Finally, the friction effect analysis yields results of mere existence, there are no characterizations or constructive clues to find the minimizing set $\omega_{friction}$. 

\item 
The difficulties presented in the analysis by friction come mostly from the non-quadratic structure of the energy dissipation functional $\F$ defined in \eqref{Eq friction dissipation functional}; this makes it highly impractical for applications. In general, whichever considered factor other than the ones presented here, leading to an energy functional whose minimizing set can not be characterized (unlike quadratic or strictly convex functionals) is likely to be impractical for computational purposes. 

\item
In the aforementioned cases, proving the existence of the associated preferential fluid directions space plays the role of theoretical support for the modeling. However, it may be wiser to use experimental sources of information and incorporate them in the global scheme with the procedure presented in section \ref{Sec Global Space of Preferential Fluid Flow Directions} for computational applications.
\end{enumerate}
So far the analysis has been made for a porous medium with one single crack embedded. However, in the future these will be applied to a fissured system, defining a Markov Chain in the following way:
\begin{enumerate}[(i)]
\item
Define a grid on the system so that each portion contains at most one crack; see figure \ref{Fig Transmission Probabilities} (a).

\item
On each element the space of preferential fluid flow directions can be defined using the procedure of this work. 

\item
Let $N$ be the number of elements in the grid. For each element, namely $k$, we define its transmission probabilities $\{p_{k, n}: 1\leq k\leq N\}$ as follows.  We set $0$ as probability transmission between two elements which do not share a common face. For each face of contact (4 faces on 2-D and 6 faces on 3-D, at most) the transmission probability from the element $k$ to its neighbor is given by the probability that the common face has to be hit by the stream line of a preferential flow direction, starting from the centroid of the element $k$; see \ref{Fig Transmission Probabilities} (b). 

\item
Since the transmission probabilities add up to one and they are all non-negative for each element, the procedure described above clearly defines a stochastic matrix. The matrix has null entries on the diagonal. It is also sparse because it has at most 4 non-null entries in 2-D and at most 6 non-null entries in 3-D.   

\end{enumerate}
%
%
\section*{Acknowledgments}
%
%
This work was supported by Universidad Nacional de Colombia.



%


\end{document}